\documentclass[11pt]{article}

\usepackage{amsmath}
\usepackage{amsfonts}
\usepackage{color}
\usepackage{hyperref}
\usepackage{multirow}
\usepackage{graphicx}
\usepackage{enumitem}
\usepackage{array}
\usepackage{calc}
\usepackage{fnpct}

\newtheorem{teo}{Theorem}[section]
\newtheorem{lem}{Lemma}[section]

\newtheorem{assump}{Assumption}

\newenvironment{pro}{\noindent\textit{Proof:}}{\halmos}
 
\DeclareMathOperator*{\Minimize}{Minimize}

\newcommand{\R}{\mathbb{R}}
\newcommand{\N}{\mathbb{N}}   
\newcommand{\nblocks}{n_{\mathrm{blocks}}}
\newcommand{\nops}{n_{\mathrm{ops}}}     
\newcommand{\nge}{n_{\mathrm{g}}}     
\newcommand{\trial}{{\mathrm{trial}}} 
\newcommand{\bound}{{\mathrm{bound}}} 
\newcommand{\chosen}{{\mathrm{chosen}}} 
\newcommand{\bx}{\boldsymbol{x}}
\newcommand{\halmos}{\hfill$\Box$}
\newcommand{\half}{\frac{1}{2}}

\begin{document}

\title{Block Coordinate Descent for smooth nonconvex constrained
  minimization\thanks{This work was supported by FAPESP (grants
    2013/07375-0, 2016/01860-1, and 2018/24293-0) and CNPq (grants
    302538/2019-4 and 302682/2019-8).}\textsuperscript{ ,}\thanks{This
    work was presented by J. M. Mart\'{\i}nez at the 5th China-Brazil
    Symposium on Applied and Computational Mathematics, that was held
    in Songshan Lake Science City from August 23rd to 24th, 2021.}}

\author{
  E. G. Birgin\thanks{Department of Computer Science, Institute of
    Mathematics and Statistics, University of S\~ao Paulo, Rua do
    Mat\~ao, 1010, Cidade Universit\'aria, 05508-090, S\~ao Paulo, SP,
    Brazil. e-mail: egbirgin@ime.usp.br}
  \and
  J. M. Mart\'{\i}nez\thanks{Department of Applied Mathematics,
    Institute of Mathematics, Statistics, and Scientific Computing
    (IMECC), State University of Campinas, 13083-859 Campinas SP,
    Brazil. e-mail: martinez@ime.unicamp.br}}

\date{November 13, 2021}

\maketitle

\begin{abstract}
At each iteration of a Block Coordinate Descent method one minimizes
an approximation of the objective function with respect to a generally
small set of variables subject to constraints in which these variables
are involved. The unconstrained case and the case in which the
constraints are simple were analyzed in the recent literature. In this
paper we address the problem in which block constraints are not simple
and, moreover, the case in which they are not defined by global sets
of equations and inequations. A general algorithm that minimizes
quadratic models with quadratric regularization over blocks of
variables is defined and convergence and complexity are proved. In
particular, given tolerances $\delta>0$ and $\varepsilon>0$ for
feasibility/complementarity and optimality, respectively, it is shown
that a measure of $(\delta,0)$-criticality tends to zero; and the the
number of iterations and functional evaluations required to achieve
$(\delta,\varepsilon)$-criticality is $O(\varepsilon^2)$. Numerical
experiments in which the proposed method is used to solve a continuous
version of the traveling salesman problem are presented.\\
  

\end{abstract}

\section{Introduction} \label{introduction}

The structure of many practical problems suggests the employment of
Block Coordinate Descent (BCD) methods for Optimization. At each
iteration of a BCD method only a block of variables is modified with
the purpose of obtaining sufficient decrease of the objective
function.

Wright~\cite{wright} surveyed traditional approaches and modern
advances on the definition and analysis of Coordinate Descent
methods. His analysis addresses mostly unconstrained problems in which
the objective function is convex. Although the Coordinate Descent
paradigm is very natural and is implicitly used in different
mathematical contexts, a classical example by Powell \cite{powell}
showed that convergence results cannot be achieved under excessively
naive implementations.

In a recent report \cite{aabmm} it was shown that, requiring
sufficient descent based on regularization, the drawbacks represented
by Powell's example can be removed. In that paper it was also shown
that methods based on high-order regularization can be defined in
which convergence and worst-case complexity can be proved. However,
the main results shown in \cite{aabmm} indicate that it is not
worthwhile to use Taylor-like models of order greater than 2 because
complexity is dominated by the necessity of keeping consecutive
iterations close enough, a requirement that is hard to achieve if
models and regularizations are of high order.  This is the reason why,
in the present paper, we restrict ourselves to quadratic models of the
objective function and quadratic regularization.

The novelty of our approach relies in the employment of a general
feasible set for each block of variables. As a consequence, at each
iteration of BCD we minimize a problem with (probably) a small number
of variables that must satisfy arbitrary constraints. Moreover, the
block feasible set may not be defined by a global set of equalities
and inequalities, as usually in constrained optimization.  Instead,
equalities and inequalities that define the feasible set are local in
nature in a sense that will be defined below, making it possible more
general domains than the ones defined by global equalities and
inequalities.
  
This paper is organized as follows. In Section~\ref{problem} the
definition of the optimization problem is given. In Section
\ref{bcsmethod} we define the BCD method for solving the main
problem. In Section~\ref{convergence} we prove convergence and
complexity results.  In Section~\ref{subpros2} we explain how to solve
subproblems. Experiments are shown in Section~\ref{experiments} and in
Section~\ref{conclusions} we state conclusions and lines for future
research.\\

\noindent
\textbf{Notation.} $\|\cdot\|$ denotes the Euclidean norm. $\nabla_i$
denotes the gradient with respect to the $i$th block of
coordinates. $\N = \{0, 1, 2,\ldots\}$.

\section{The problem} \label{problem}

Assume that $f: \R^n \to \R$, $n_i \geq 1$ for all
$i=1,\ldots,\nblocks$, and $\sum_{i=1}^{\nblocks} n_i = n$.  Let us
write $x^T = (\bx_1^T, \dots, \bx_{\nblocks}^T)^T$, where $\bx_i \in
\R^{n_i}$ for all $i=1, \dots, \nblocks$, so
\[
f(x) = f(x_1, \dots, x_n) = f(\bx_1, \dots, \bx_{\nblocks}).
\]
The problem considered in the present work is given by
\begin{equation} \label{theproblem}
\Minimize f(\bx_1, \dots, \bx_{\nblocks}) \mbox{ subject to }
\bx_i \in \Omega_i, i = 1,\dots,\nblocks.
\end{equation}

The assumption below includes conditions on the sets $\Omega _i$ and
on the way they are described.

\begin{assump} \label{a0}    
For all $i=1,\dots,\nblocks$ the set $\Omega_i \subset \R^{n_i}$ is
closed and bounded. Moreover, for all $i = 1, \dots,\nblocks$, there
exist open sets $A_{i,j}$, $j=1,\dots,\nops(i)$, such that
\[
\Omega_i \subset \cup_{j=1}^{\nops(i)} A_{i,j}
\]
and there exist $\nge(i,j)$ smooth functions
\begin{equation} \label{ineq1}
g_{i,j,\ell}: A_{i,j} \to \R, \; \ell=1,\dots,\nge(i,j)
\end{equation}
such that  
\begin{equation} \label{ineq2}
\Omega_i \cap A_{i,j} = \{\bx_i \in A_{i,j} \;|\; g_{i,j,\ell}(\bx_i)
\leq 0, \; \ell=1,\dots,\nge(i,j)\}.
\end{equation}

\end{assump}

\begin{figure}[ht!]
\begin{center}
\includegraphics{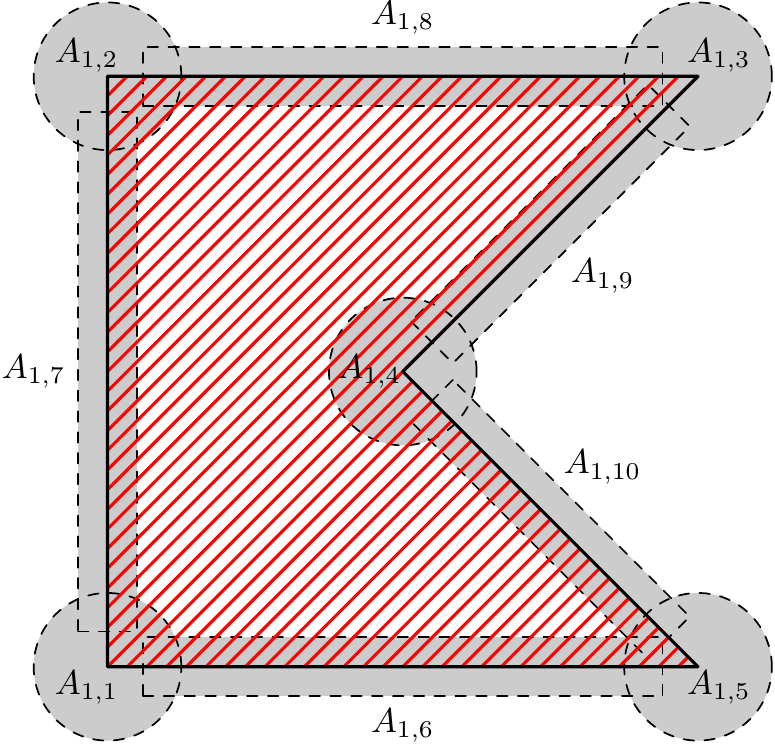}
\end{center}
\caption{Illustration of a set $\Omega_1$ (dashed red) covered by
  $\nops(1)=12$ open sets $A_{1,1}, \dots, A_{1,12}$.}
\label{fig1}
\end{figure}

The constraints $g_{i,j,\ell}(\bx_i) \leq 0$ are said to be
\textit{constitutive} of $\Omega_i$ in the open covering set
$A_{i,j}$. The explicit inclusion of equality constraints
in~(\ref{ineq1},\ref{ineq2}) offers no difficulty and we state the
case with only inequalities for the sole purpose of simplifying the
notation. See Figure~\ref{fig1} for an example of a set $\Omega_i$,
the open covering $A_{i,j}$ and functions $g_{i,j,\ell}$.

In Figure~\ref{fig1}, sets $A_{1,j}$ for $j=1,2,3,5$ (four out of the
five open circles) are such that $n_g(1,j)=2$; while sets $A_{1,j}$
for $j=6,7,8,9,10$ (the open rectangles) are such that
$n_g(1,j)=1$. In all cases, constitutive constraints are linear. Sets
$A_{1,11}$ and $A_{1,12}$ are \textit{not} displayed in the picture
for clarity. They are two congruent open triangles that cover the
interior of $\Omega_1$ that appears unconvered in the picture; and
they are such that $n_g(1,11)=n_g(1,12)=0$. The ``internal kink'' at
$A_{1,4}$ makes the constitutive constraint of $A_{1,4}$ to deserve
special consideration.  In $A_{1,4}$ the feasible region is the
complement of a set defined by constraints of the form $a_1 x_1 + b_1
x_2 + c_1 \geq 0$ and $a_2 x_1 + b_2 x_2 + c_2 \geq 0$. Let us write
$z_1 = a_1 x_1 + b_1 x_2 + c_1$ and $z_2 = a_2 x_1 + b_2 x_2 +
c_2$. Then, in the plane $(z_1, z_2)$ the feasible region is, locally,
the complement of the non-negative orthant. Define $\varphi(z_1, z_2)
= (z_1 z_2)^2$ if $z_1 \geq 0$ and $z_2 \geq 0$, whereas $\varphi(z_1,
z_2) = - (z_1 z_2)^2$ otherwise.  Then, the feasible region is,
locally, given by $\varphi(z_1, z_2) \leq 0$; i.e.\ $n_g(1,4)=1$ and
the constitutive constraint is given by~$\varphi$.

For further reference (in Section~\ref{subpros2}), we show here that
the center $C_{1,4}$ of the ball $A_{1,4}$ satisfies KKT
conditions. Clearly, the origin in the plane $(z_1,z_2)$, which
corresponds to the point $C_{1,4}$, is a non-regular point from the
point of view of constrained optimization. However, if a smooth
function $\psi$ has a minimizer at this point, its gradient $\nabla
\psi(A_{1,4})$ is necessarily null. To see this observe that it is
easy to show two linearly independent directions $v_1$ and $v_2$ such
that $\nabla \psi(C_{1,4})^T v_1 \geq 0$, $\nabla \psi(C_{1,4})^T
(-v_1) \geq 0$, $\nabla \psi(C_{1,4})^T v_2 \geq 0$, and $\nabla
\psi(C_{1,4})^T (-v_2) \geq 0$, implying that $\nabla \psi(C_{1,4}) =
0$, so that $C_{1,4}$ is a KKT point of the minimization of $\psi$
subject to $\varphi \leq 0$ with a single null multiplier.

\section{Block Coordinate Descent method} \label{bcsmethod}

In this section we present the main algorithm designed for solving
(\ref{theproblem}). At each iteration $k$ of the Block Coordinate
Descent (BCD) method, given the current iterate $\bx^k =((\bx_1^k)^T,
\dots, (\bx_{\nblocks}^k)^T)^T$, we choose $i_k \in \{1, \dots,
\nblocks\}$ and compute $\bx_{i_k}^{k+1}$ by approximately solving
\begin{equation} \label{solvapro}
\Minimize_{\bx_{i_k} \in \R^{n_i}} f(\bx_1, \dots, \bx_{\nblocks})
\mbox{ subject to } \bx_{i_k} \in \Omega_{i_k} \mbox{ and } \bx_j =
\bx_j^k \mbox{ for all } j \neq i_k,
\end{equation}
i.e.\ we approximately minimize the function $f$ fixing the blocks
$\bx_j$ such that $j \neq i_k$. For $j \neq i_k$, we define
$\bx_j^{k+1} = \bx_j^k$. The sense in which problem (\ref{solvapro})
is solved only approximately is clarified below.

The algorithmic parameters of BCD are the sufficient descent parameter
$\alpha$, which defines progress of the objective function and
implicitly penalizes the distance between consecutive iterates, the
tolerance~$\delta$ with respect to complementarity conditions, the
parameter~$\theta$ that defines sufficient descent of the model at
each iteration, and the minimal positive regularization
parameter~$\sigma_{\min}$. The model Hessian matrices~$B_k$ may be
used to mimic available second derivative information but do not play
any significative role from the point of view of complexity or
convergence and the choice $B_k = 0$ is always possible. At Step~2 of
the algorithm we choose the block of variables with respect to which
we wish to improve the objective function. In general, we minimize
approximately a quadratic model increasing the regularizing parameter
as far as the suffcient condition (\ref{armijo}) is not
satisfied. Alternatively, we employ the true objective function as a
model, because such alternative is possible in many practical
problems.\\

\noindent
\textbf{Algorithm \ref{bcsmethod}.1.} Assume that $\alpha > 0$,
$\delta \geq 0$, $\theta > 0$, $\sigma_{\min} > 0$, and $\bx_i^0 \in
\Omega_i$ for $i=1,\dots,\nblocks$ are given. Initialize the iteration
number $k \leftarrow 0$.
\begin{description}
\item[\textbf{Step~1.}] Set $\sigma \leftarrow 0$, choose
  $i_k \in \{1,\dots,\nblocks\}$, and define $B_k \in \R^{n_{i_k}}
  \times \R^{n_{i_k}}$ symmetric.
  
\item[\textbf{Step~2.}] Find $j_k \in \{1,\dots,\nops(i_k)\}$ and
  $\bx_{i_k}^{\trial} \in A_{i_k,j_k}$ such that if $\sigma>0$ then
  Alternative~1 holds while if $\sigma=0$ then either Alternative~1 or
  Alternative~2 holds.

\textbf{Alternative 1:}
\begin{equation} \label{bajar} 
  \nabla_{i_k} f(x^k)^T (\bx_{i_k}^{\trial}-\bx_{i_k}^k) + \half
  (\bx_{i_k}^{\trial}-\bx_{i_k}^k)^T B_k
  (\bx_{i_k}^{\trial}-\bx_{i_k}^k) + \frac{\sigma}{2} \|\bx_{i_k}^{\trial}
  - \bx_{i_k}^k\|^2 \leq 0
\end{equation}            
and   there exist $\mu_{i_k,j_k,\ell} \geq
0$ for $\ell = 1,\dots, \nge(i_k,j_k)$ for which
\begin{equation} \label{grachico}
  \left\| \nabla_{i_k} f(x^k) + B_k (\bx_{i_k}^{\trial}-\bx_{i_k}^k)
  + \sigma (\bx_{i_k}^{\trial} - \bx_{i_k}^k) +
  \sum_{\ell=1}^{\nge(i_k,j_k)} \mu_{i_k,j_k,\ell} \nabla
  g_{i_k,j_k,\ell}(\bx_{i_k}^{\trial}) \right\| \leq \theta
  \|\bx_{i_k}^{\trial}-\bx_{i_k}^k\|
\end{equation}
and      
\begin{equation} \label{lambdamu}
  \min\{\mu_{i_k,j_k,\ell}, -g_{i_k,j_k,\ell}(\bx_{i_k}^{\trial})\} \leq \delta, \;
  \ell=1,\dots,\nge(i_k, j_k).
\end{equation}

\textbf{Alternative 2:}

There exist $\mu_{i_k,j_k,\ell} \geq 0$ for $\ell = 1,\dots,
\nge(i_k,j_k)$ for which
\begin{equation} \label{grachi2}
  \left\| \nabla_{i_k} f(\bx_1^k, \dots, \bx_{i_k}^{\trial}, \dots,
  \bx_{\nblocks}^k) + \sum_{\ell=1}^{\nge(i_k,j_k)}
  \mu_{i_k,j_k,\ell} \nabla g_{i_k,j_k,\ell}(\bx_{i_k}^{\trial})
  \right\| = 0
\end{equation}
and 
\begin{equation} \label{lamb2}
  \min\{\mu_{i_k,j_k,\ell}, -g_{i_k,j_k,\ell}(\bx_{i_k}^{\trial})\}
  \leq 0, \; \ell=1,\dots,\nge(i_k, j_k).
\end{equation}

\item[\textbf{Step~3.}] Test the sufficient descent condition
  \begin{equation} \label{armijo}   
   f(\bx_1^k, \dots, \bx_{i_k}^{\trial}, \dots, \bx_{\nblocks}^k) \leq
   f(x^k) - \alpha \|\bx_{i_k}^{\trial} - \bx_{i_k}^k\|^2.
  \end{equation}
 If (\ref{armijo}) holds, then define $\sigma_k = \sigma$,
 $x_{i_k}^{k+1} = x_{i_k}^{\trial}$ and $x_i^{k+1} = x_i^k$ for all $i
 \neq i_k$, set $k \leftarrow k+1$ and go to Step~1. Otherwise, set
 $\sigma \leftarrow \max\{ \sigma_{\min}, 2 \sigma \}$ and go to
 Step~2.
\end{description}

\noindent
\textbf{Remark.} We will see that, from the theoretical point of view,
Alternative 2 is not necessary. Convergence and complexity theoretical
results follow without any difficulty with Alternative~1
only. Alternative~2 was included because, in many cases, a procedure
exists to find a global minimizer with respect to a single block. So,
in Alternative 2 we allow the algorithm to choose such minimizer, with
the only condition that it must satisfy KKT conditions in the
block. However, note that the test~(\ref{armijo}) is still necessary
and cannot be eliminated. The reason is that its fulfillment implies
that the difference between consecutive iterations tends to zero and
this feature is essential for the convergence of coordinate search
methods. See the counterexample in~\cite{powell} and the discussion
with only box constraints in~\cite{aabmm}.
  
\section{Convergence and complexity} \label{convergence}

In this section we prove convergence and complexity results. We say
that $(\bx_1, \dots, \bx_{\nblocks})$ is
$(\delta,\varepsilon)$-critical if there exist open sets $A_{i,j}$
($i=1,\dots,\nblocks$, $j=1,\dots,\nops(i)$) satisfying
Assumption~\ref{a0} such that, for all $i=1,\dots,\nblocks$, $\bx_i$
satisfies the KKT conditions for the minimization of
$f(\bx_1,\ldots,\bx_{\nblocks})$ restricted to the constitutive
constraints of $A_{i,j}$, $j=1\dots,\nops(i)$, with tolerance
$\varepsilon>0$ and satisfies complementarity and feasibility with
respect to the same constraints with tolerance $\delta>0$. Under
proper assumptions we prove that, given $\delta \geq 0$ (which is a
parameter of the algorithm BCD), the natural measure of $(\delta,
0)$-criticality tends to zero and the number of iterations and
evaluations that are necessary to obtain $(\delta,
\varepsilon)$-criticality is $O(\varepsilon^2)$.

In Assumption~\ref{a1} we state that the gradients of the objective
function must satisfy Lipschitz conditions.

\begin{assump} \label{a1}
There exists $\gamma > 0$ such that, for all $\bx_{i_k}$ and
$\bx_{i_k}^{\trial}$ computed at Algorithm~\ref{bcsmethod}.1,
\begin{equation} \label{lipschitz}
\| \nabla_{i_k} f(x^k) - \nabla_{i_k} f(\bx_1^k, \dots,
\bx_{i_k}^{\trial}, \dots, \bx_{\nblocks}^k) \| \leq \gamma
\|\bx_{i_k} - \bx_{i_k}^{\trial}\|
\end{equation}
and
\begin{equation} \label{taylor}      
f(\bx_1^k, \dots, \bx_{i_k}^{\trial}, \dots, \bx_{\nblocks}^k) \leq
f(x^k) + \nabla_{i_k} f(x^k)^T (\bx_{i_k}^{\trial}-\bx_{i_k}^k) +
\frac{\gamma}{2} \|\bx_{i_k} - \bx_{i_k}^{\trial}\|^2.
\end{equation}
Moreover, if $x^{k_1}$ differs from $x^{k_2}$ in only one block
$i_{\mathrm{diff}}$,
\begin{equation} \label{lips2}
\|\nabla_{i} f(x^{k_1}) - \nabla_{i} f(x^{k_2})\| \leq \gamma
\|\bx_{i_{\mathrm{diff}}}^{k_2} - \bx_{i_{\mathrm{diff}}}^{k_1}\|.
\end{equation}                                         
for all $i=1,\dots,\nblocks$.
\end{assump}

Assumption~\ref{a2} merely states that model Hessians should be
uniformly bounded.

\begin{assump} \label{a2}
There exist $c_B > 0$ such that for all $ k \in \N$,
\begin{equation} \label{cotab}
\|B_k\| \leq c_B.
\end{equation}
\end{assump}

Assumptions~\ref{a1} and~\ref{a2} are sufficient to prove that every
iteration of BCD is well defined, as sufficient descent (\ref{armijo})
is obtained increasing the regularization parameter $\sigma$ a finite
number of times.

\begin{lem} \label{lem1}
Assume that Assumptions~\ref{a1} and~\ref{a2} hold and that, for all
$k \in \N$, the computation of $j_k$ and $\bx_{i_k}^{\trial}$
according to Step~2 of BCD is possible. Then, the test (\ref{armijo})
is satisfied after at most
\[
\log_2 \left( \frac{\gamma + c_B + 2 \alpha}{\sigma_{\min}} \right) + 1
\]
increases of $\sigma$ at Step~3. Moreover,
\begin{equation} \label{cotasigma}
\sigma_k < \sigma_{\max} := 2 (\gamma + c_B + 2 \alpha).
\end{equation}
\end{lem}

\begin{pro}
If the test~(\ref{armijo}) is satisfied when Step~2 is executed with
$\sigma=0$, then the thesis holds trivially. So, we need to consider
only the case in which $\sigma>0$ and, in consequence, $j_k$
and~$\bx_{i_k}^{\trial}$ satisfy Alternative~1.
By (\ref{taylor}) in Assumption~\ref{a1},
\begin{align*}
f(\bx_1^k, \dots, \bx_{i_k}^{\trial}, \dots, \bx_{\nblocks}^k) &\leq
f(x^k) + \nabla_{i_k} f(x^k)^T (\bx_{i_k}^{\trial}-\bx_{i_k}^k) +
\frac{\gamma}{2} \|\bx_{i_k} - \bx_{i_k}^{\trial}\|^2\\
&+ \left( \half (\bx_{i_k}^{\trial}-\bx_{i_k}^k)^T B_k
(\bx_{i_k}^{\trial}-\bx_{i_k}^k) + \half \sigma \|\bx_{i_k}^{\trial} -
\bx_{i_k}^k\|^2 \right)\\
&- \left( \half (\bx_{i_k}^{\trial}-\bx_{i_k}^k)^T B_k
(\bx_{i_k}^{\trial}-\bx_{i_k}^k) + \half \sigma \|\bx_{i_k}^{\trial} -
\bx_{i_k}^k\|^2 \right)
\end{align*}
Then, by (\ref{bajar}),
\[
f(\bx_1^k, \dots, \bx_{i_k}^{\trial}, \dots, \bx_{\nblocks}^k) \leq
f(x^k) + \frac{\gamma}{2} \|\bx_{i_k} - \bx_{i_k}^{\trial}\|^2 - \half
(\bx_{i_k}^{\trial}-\bx_{i_k}^k)^T B_k
(\bx_{i_k}^{\trial}-\bx_{i_k}^k) - \half \sigma \|\bx_{i_k}^{\trial} -
\bx_{i_k}^k\|^2.
\]
Therefore, by (\ref{cotab}) in Assumption~\ref{a2}, 
\begin{align*}
f(\bx_1^k, \dots, \bx_{i_k}^{\trial}, \dots, \bx_{\nblocks}^k) &\leq
f(x^k) + \frac{\gamma}{2} \|\bx_{i_k} - \bx_{i_k}^{\trial}\|^2 + \half
c_B \| \bx_{i_k}^{\trial}-\bx_{i_k}^k \|^2 - \half \sigma
\|\bx_{i_k}^{\trial} - \bx_{i_k}^k\|^2\\
&= f(x^k) + \half ( \gamma + c_B - \sigma ) \|\bx_{i_k}^{\trial} - \bx_{i_k}^k\|^2.
\end{align*}
Then, the inequality (\ref{armijo}) holds if
\[
\half ( \gamma + c_B - \sigma ) \|\bx_{i_k}^{\trial} - \bx_{i_k}^k\|^2 \leq
-\alpha \|\bx_{i_k}^{\trial} - \bx_{i_k}^k\|^2,
\]
i.e.\ if $\sigma \geq \gamma + c_B + 2 \alpha$. Since, by definition,
$\sigma$ initially receives the value zero and then receives values of
the form $2^{\ell-1} \sigma_{\min}$, where $\ell$ is the number of
executions of $\sigma \leftarrow \max\{\sigma_{\min}, 2 \sigma\}$,
then the number of increases of $\sigma$ that are necessary to
obtain~(\ref{armijo}) is bounded above by
\[
\log_2 \left( \frac{\gamma + c_B + 2 \alpha}{\sigma_{\min}} \right) + 1
\]
as we wanted to prove. Finally, (\ref{cotasigma}) comes from the fact
that the largest unsuccessful value of~$\sigma$ must be strictly less
than $\gamma + c_B + 2 \alpha$ and the next (successful) value is
twice that amount by definition.
\end{pro}

\begin{assump} \label{a3}
The sequence $\{f(\bx^k)\}$ is bounded below.
\end{assump}

The following lemma is a simple consequence of (\ref{armijo}) and
Assumption~\ref{a3}.
   
\begin{lem} \label{lem2}
Assume that Assumptions~\ref{a1}, \ref{a2}, and~\ref{a3}
hold and that for all $k \in \N$, the computation of $j_k$ and
$\bx_{i_k}^{\trial}$ at Step~2 is possible. Then, $\lim_{k \to \infty}
\|\bx_{i_k}^{k+1} - \bx_{i_k}^k\| = \lim_{k\to\infty}
\|\bx^{k+1}-\bx^k\| = 0$ and, given $\varepsilon > 0$, the number of
iterations at which $\|\bx_{i_k}^{k+1} - \bx_{i_k}^k\| > \varepsilon$
is bounded above by
\begin{equation} \label{boundif} 
\frac{f(\bx^0) - f_{\bound}}{\alpha \varepsilon^2}
\end{equation}
where $f_{\bound}$ is an arbitrary lower bound of $\{f(\bx^k)\}$.
\end{lem}

\begin{pro}
By Assumption~\ref{a3} there exists $f_{\bound} \in \R$ such that
$f(\bx^k) \geq f_{\bound}$ for all $k \in \N$. Then, the fact that
$\|\bx^{k+1}_{i_k}-\bx^k_{i_k}\|$ tends to zero comes
from~(\ref{armijo}); and this implies that $\|\bx^{k+1}-\bx^k\|$ tends
to zero as well because, by definition, $\|\bx^{k+1}-\bx^k\| =
\|\bx^{k+1}_{i_k}-\bx^k_{i_k}\|$. Finally, if $\|\bx_{i_k}^{k+1} -
\bx_{i_k}^k\| > \varepsilon$, then, by~(\ref{armijo}), $f(x^{k+1})
\leq f(x^k) - \alpha \varepsilon^2$; and this reduction can no occur
more than~(\ref{boundif}) times, as we wanted to prove.
\end{pro}\\

Assumption~\ref{a4} states that every block of components $i$ is
chosen for minimization at infinitely many iterations and, moreover,
at every $m$ consecutive iterations we necessarily find at least one
at which $i$ is chosen.
 
\begin{assump} \label{a4}
There exists $m \in \{1,2,\dots\}$ such that, for all $\nu \in
\{1,\dots,\nblocks\}$, $i_k = \nu$ infinitely many times and, if $k_1
< k_2 < k_3, \dots$ is the set of all the iteration indices $k$ such
that $i_k=\nu$, one has that $k_1 \leq m$, and $k_{j+1}-k_j \leq m$
for all $j= 1, 2, 3, \dots$.
\end{assump}

The following theorems are the main convergence result of this
paper. The idea is the following. According to
Algorithm~\ref{bcsmethod}.1, at iteration~$k$, we select a block $i_k
= i_{\chosen}$ and optimize with respect to the variables of this
block up to the approximate fulfillment of restricted KKT
conditions. These restricted KKT conditions hold in one of the open
sets that cover $\Omega_{i_{\chosen}}$ and involve the constraints
that are constitutive in this open set. The variables
$\bx_{i_{\chosen}}$ do not change during some (less than $m$)
iterations; therefore, during these iterations, thanks to the
Lipschitz assumption~(\ref{lips2}), the approximate KKT conditions
with respect to the variables $i_{\chosen}$ still hold with respect to
the same open set and the same constitutive constraints used at
iteration~$k$.  After these (less than $m$) iterations the block
$i_{\chosen}$ is selected again, and the process is repeated.  Since
all the blocks are chosen infinitely many times in the way described
by Assumption~\ref{a4}, approximate KKT conditions eventually hold
with respect to all the blocks and we are able to establish the number of
iterations that we need for the fulfillment of KKT conditions up to an
arbitrary precision~$\varepsilon$.
  
\begin{teo} \label{teo1}   
Assume that Assumptions~\ref{a1}, \ref{a2}, \ref{a3}, and~\ref{a4}
hold and that for all $k \in \N$, the computation of $j_k$ and
$\bx_{i_k}^{\trial}$ at Step~2 is possible. For $i \in \{1, \dots,
\nblocks\}$ and $k \geq m$, define $o(i,k):=j_{\hat k}$, where $\hat k$
is the latest iteration (not larger than $k$) at which
$i_{\hat k}=i$. Then, for $i \in \{1, \dots, \nblocks\}$ and $k
\geq m$, we have that $\mu_{i,o(i,k),\ell} \geq 0$ for $\ell =
1,\dots,\nge(i,o(i,k))$,
\[
\bx^k_i \in A_{i,o(i,k)},
\]
\[
\min \{\mu_{i,o(i,k),\ell}, -g_{i,o(i,k),\ell}(\bx^k_i) \} \leq \delta,
\]
and
\[
\lim_{k\to\infty} \left\| \nabla_i f(x^k) +
\sum_{\ell=1}^{\nge(i,o(i,k))} \mu_{i,o(i,k),\ell} \nabla
g_{i,o(i,k),\ell}(\bx^k_i) \right\| = 0.
\]
Moreover, given $\varepsilon > 0$,
the number of iterations at which
\[
\left\| \nabla_i f(x^k) + \sum_{\ell=1}^{\nge(i,o(i,k))}
\mu_{i,o(i,k),\ell} \nabla g_{i,o(i,k),\ell}(\bx^k_i) \right\| >
\varepsilon
\]
is bounded above by
\[
\frac{f(\bx^0) - f_{\bound}}{\alpha (\varepsilon/(c_4 m))^2},
\]
where $f_{\bound}$ is an arbitrary lower bound of $\{f(\bx^k)\}$ and
\begin{equation} \label{c4}
c_4 := c_B + \sigma_{\max} + \theta + \gamma.
\end{equation}
\end{teo}    

\begin{pro}
Let $i \in \{1, \dots, \nblocks\}$ be arbitrary. By
Assumption~\ref{a4}, there exists $k_1 \leq m$ such that $i_{k_1} =
i$. Without loss of generality, in order to simplify the notation,
suppose that $k_1 = 0$. Consider first the case where, in Step~2,
Alternative~1 holds. Then, at iteration $k_1 = 0$ one defines $B_0$
and finds $j_0$, $\bx_i^1 \in A_{i,j_0}$, and $\mu_{i,j_0,\ell} \geq
0$ for $\ell = 1,\dots, \nge(i,j_0)$ such that
\begin{equation} \label{grachi}
\left\| \nabla_i f(x^0) + B_0 (\bx_i^1-\bx_i^0) + \sigma_0 (\bx_i^1 -
\bx_i^0) + \sum_{\ell=1}^{\nge(i,j_0)} \mu_{i,j_0,\ell} \nabla
g_{i,j_0,\ell}(\bx_i^1) \right\| \leq \theta \|\bx_i^1-\bx_i^0\|
\end{equation}
and
\begin{equation} \label{lamb}
\min\{\mu_{i,j_0,\ell}, -g_{i,j_0,\ell}(\bx_i^1)\} \leq \delta
\;\mbox{ for }\; \ell=1,\dots, \nge(i,j_0).
\end{equation}
By (\ref{grachi}) and (\ref{cotasigma}),  
\begin{equation} \label{grach}
\left\| \nabla_i f(x^0) + \sum_{\ell=1}^{\nge(i,j_0)} \mu_{i,j_0,\ell}
\nabla g_{i,j_0,\ell}(\bx_i^1) \right\| \leq (c_B + \sigma_{\max}
+\theta) \|\bx_i^1-\bx_i^0\|.
\end{equation}
Then, by (\ref{lipschitz}),
\begin{equation} \label{grac}
\left\| \nabla_i f(\bx_1^0, \dots, \bx_i^1, \dots, \bx_{\nblocks}^0) +
\sum_{\ell=1}^{\nge(i,j_0)} \mu_{i,j_0,\ell} \nabla
g_{i,j_0,\ell}(\bx_i^1) \right\| \leq (c_B + \sigma_{\max} +\theta +
\gamma) \|\bx_i^1-\bx_i^0\|.
\end{equation}  
Since $\bx_s^1 = \bx_s^0$ for every $s=1,\dots,\nblocks$, $s \neq i$,
by (\ref{grac}) and the definition of $c_4$ in (\ref{c4}), we have
that
\begin{equation} \label{gr}
\left\| \nabla_i f(x^1) + \sum_{\ell=1}^{\nge(i,j_0)}
\mu_{i,j_0,\ell} \nabla g_{i,j_0,\ell}(\bx_i^1) \right\| \leq c_4
\|\bx_i^1-\bx_i^0\|.
\end{equation}       
Recall that~(\ref{lamb}) and~(\ref{gr}) were obtained under
Alternative~1. On the other hand, under Alternative~2, (\ref{lamb})
and~(\ref{gr}) follow trivially from~(\ref{lamb2})
and~(\ref{grachi2}), respectively.

Since $\bx^2$ may differ from $\bx^1$ only in the block $i_1$,
by~(\ref{lips2}), we have that
\[
\|\nabla_i f(\bx^2) - \nabla_i f(\bx^1)\| \leq \gamma \|\bx_{i_1}^2 -
\bx_{i_1}^1\|.
\]
Then, by (\ref{gr}),
\begin{equation} \label{gr2}
\left\| \nabla_i f(x^2) + \sum_{\ell=1}^{\nge(i,j_0)}
\mu_{i,j_0,\ell} \nabla g_{i,j_0,\ell}(\bx_i^1) \right\| \leq c_4
\|\bx_i^1-\bx_i^0\| + \gamma \|\bx_{i_1}^2 - \bx_{i_1}^1\|.
\end{equation}     
Since $\bx^3$ may differ from $\bx^2$ only in the block $i_2$,
by~(\ref{lips2}) and~(\ref{gr2}), we have that
\begin{equation} \label{gr3}
\left\| \nabla_i f(x^3) + \sum_{\ell=1}^{\nge(i,j_0)}
\mu_{i,j_0,\ell} \nabla g_{i,j_0,\ell}(\bx_i^1) \right\| \leq c_4
\|\bx_i^1-\bx_i^0\| + \gamma \|\bx_{i_1}^2 - \bx_{i_1}^1\| + \gamma
\|\bx_{i_2}^3 - \bx_{i_2}^2\|.
\end{equation}     
So, using an inductive argument, for all $k \in \N$,
\begin{equation} \label{grkk}        
\left\| \nabla_i f(x^k) + \sum_{\ell=1}^{\nge(i,j_0)}
\mu_{i,j_0,\ell} \nabla g_{i,j_0,\ell}(\bx_i^1) \right\| \leq c_4
\|\bx_i^1-\bx_i^0\| + \gamma \sum_{\nu=2}^k \|\bx^\nu - \bx^{\nu-1}\|.
\end{equation}
Thus, by the definition of $c_4$ in (\ref{c4}), for all $k \in \N$,
\begin{equation} \label{grkkk}        
\left\| \nabla_i f(x^k) + \sum_{\ell=1}^{\nge(i,j_0)}
\mu_{i,j_0,\ell} \nabla g_{i,j_0,\ell}(\bx_i^1) \right\| \leq c_4
\sum_{\nu=1}^k \|\bx^\nu - \bx^{\nu-1}\|.
\end{equation} 

Let us now get rid of the simplifying assumption $i_0 = i$ and assume
that the set of indices~$k$ at which $i_k=i$ is $k_1 < k_2 <
\dots$. Renaming the indices in (\ref{lamb}) and (\ref{grkkk}), we get
that
\begin{equation} \label{lambnovo}
\min\{\mu_{i,j_{k_r},\ell}, -g_{i,j_{k_r},\ell}(\bx_i^{k_r+1})\} \leq \delta
\;\mbox{ for }\; \ell=1,\dots, \nge(i,j_{k_r})
\end{equation}
and for all $r=1,2,\dots$ and all $k=k_r+1,k_r+2,\dots$, in particular
for all $k=k_r+1,k_r+2,\dots,k_{r+1}$,
\begin{equation} \label{grok1}        
\left\| \nabla_i f(x^k) + \sum_{\ell=1}^{\nge(i,j_{k_r})}
\mu_{i,j_{k_r},\ell} \nabla g_{i,j_{k_r},\ell}(\bx_i^{k_r+1}) \right\|
\leq c_4 \sum_{\nu=k_r+1}^{k} \|\bx^\nu - \bx^{\nu-1}\|.
\end{equation}    
But, by the definition of the sequence $k_1, k_2, \dots$, $\bx_i^{k_r}$
may change from iteration $k_r$ to iteration $k_r+1$ but it does not
change from $k_r+1$ to $k_r+2$, $k_r+2$ to $k_r+3$, until $k_{r+1}-1$
to $k_{r+1}$; and it may change again from iteration $k_{r+1}$ to
iteration $k_{r+1}+1$. This means that, for all $r=1,2,\dots$, we have that 
$\bx_i^{k_r+1} = \bx_i^{k_r+2} = \dots = \bx_i^{k_{r+1}}$. Therefore,
(\ref{lambnovo}) and~(\ref{grok1}) imply that for all $r=1,2,\dots$
and $k=k_r+1, k_r+2,\dots,k_{r+1}$,
\begin{equation} \label{lambk}
\min\{\mu_{i,j_{k_r},\ell}, -g_{i,j_{k_r},\ell}(\bx_i^k)\} \leq \delta
\;\mbox{ for }\; \ell=1,\dots, \nge(i,j_{k_r})
\end{equation}
and
\begin{equation} \label{grip}        
\left\| \nabla_{i} f(x^k) + \sum_{\ell=1}^{\nge(i,j_{k_r})}
\mu_{i,j_{k_r},\ell} \nabla g_{i,j_{k_r},\ell}(\bx_i^{k})
\right\| \leq c_4 \sum_{\nu=k_r+1}^k \|\bx^\nu - \bx^{\nu-1}\|.
\end{equation}      
Moreover, by definition, $o(i,k_r+1) = o(i,k_r+2) = \dots =
o(i,k_{r+1}) = j_{k_r}$ for $r=1,2,\dots$. So, from (\ref{lambk})
and~(\ref{grip}), we get that, for $k \geq m$, $\bx_i^k \in
A_{i,o(i,k)}$ and $\mu_{i,j_0,\ell} \geq 0$ for $\ell = 1,\dots,
\nge(i,j_0)$ are such that
\begin{equation} \label{lambkbis}
\min\{\mu_{i,o(i,k),\ell}, -g_{i,o(i,k),\ell}(\bx_i^k)\} \leq \delta
\;\mbox{ for }\; \ell=1,\dots, \nge(i,o(i,k))
\end{equation}
and
\begin{equation} \label{gripbis}        
\left\| \nabla_{i} f(x^k) + \sum_{\ell=1}^{\nge(i,o(i,k))}
\mu_{i,o(i,k),\ell} \nabla g_{i,o(i,k),\ell}(\bx_i^{k})
\right\| \leq c_4 \sum_{\nu=k_r+1}^k \|\bx^\nu - \bx^{\nu-1}\|.
\end{equation}      

By Lemma~\ref{convergence}.2, given $\varepsilon > 0$, the number of
iterations at which $\|\bx^\nu - \bx^{\nu-1}\| > \varepsilon/(c_4 m$)
is bounded above by $( f(\bx^0) - f_{\bound} ) / ( \alpha
(\varepsilon/(c_4 m))^2 )$. Then, since the sum in the second member
of (\ref{grip}) involves at most $m$ terms, the number of iterations
at which the left-hand side of (\ref{grip}) is bigger than
$\varepsilon$ is bounded above by $( f(\bx^0) - f_{\bound} ) / (
\alpha (\varepsilon/(c_4 m))^2 )$. Moreover, since $\varepsilon>0$ is
arbitrary, taking limits on both sides of~(\ref{gripbis}), we have
that
\begin{equation} \label{gripbislim}        
\lim_{k\to\infty} \left\| \nabla_{i} f(x^k) + \sum_{\ell=1}^{\nge(i,o(i,k))}
\mu_{i,o(i,k),\ell} \nabla g_{i,o(i,k),\ell}(\bx_i^{k})
\right\| = 0.
\end{equation}      
\end{pro}

The final theorem of this section proves worst-case functional
complexity of order $O(\varepsilon^{-2})$.

\begin{teo} \label{teo2}   
Assume that Assumptions~\ref{a1}, \ref{a2}, \ref{a3}, and~\ref{a4}
hold and that for all $k \in \N$, the computation of $j_k$ and
$\bx_{i_k}^{\trial}$ at Step~2 is possible. Let the indices $o(i,k)$
be as defined in Theorem~\ref{teo1}. Then, given $\varepsilon>0$,
Algorithm~\ref{bcsmethod}.1 performs at most
\[
\nblocks \left( \frac{f(x^0) - f_{\bound}}{\alpha (\varepsilon/(c_4 m))^2} \right)
\]
iterations and at most
\[
\log_2 \left( \frac{\gamma + c_B + 2 \alpha}{\sigma_{\min}} \right) + 2
\]
functional evaluations per iteration, where $f_{\bound}$ is an
arbitrary lower bound of $\{f(\bx^k)\}$ and $c_4$ is given
by~(\ref{c4}), to compute an iterate $x^{k+1}$ such that
\[
\bx^{k+1}_i \in A_{i,o(i,k)},
\]
\[
\mu_{i,o(i,k),\ell} \geq 0 \;\mbox{ and }\; \min
\{\mu_{i,o(i,k),\ell}, -g_{i,o(i,k),\ell}(\bx^{k+1}_i) \} \leq \delta,
\; \ell = 1,\dots,\nge(i,o(i,k)),
\]
and
\[
\left\| \nabla_i f(x^{k+1}) + \sum_{\ell=1}^{\nge(i,o(i,k))}
\mu_{i,o(i,k),\ell} \nabla g_{i,o(i,k),\ell}(\bx^{k+1}_i) \right\|
\leq \varepsilon
\]
for all $i=1,\dots,\nblocks$.
\end{teo}    

\begin{pro}
The proof follows from Theorem~\ref{teo1}, Lemma~\ref{lem1}, and the
definition of Algorithm~\ref{bcsmethod}.1, because the number of
functional evaluations per iterations is equal to the number of
increments of $\sigma$ plus one. (This ignores the fact that,
disregarding the first iteration, the value of $f(x^k)$ can in fact be
obtained from the previous iteration, in which case the number of
functional evaluations and the number of increases of $\sigma$ per
iteration coincide.)
\end{pro}

\section{Solving subproblems} \label{subpros2}

In this section we present an algorithm for computing
$\bx_{i_k}^{\trial}$ at Step~2 of iteration~$k$ of
Algorithm~\ref{bcsmethod}.1. The well-definiteness of the proposed
approach requires, in addition to~(\ref{ineq2}), two additional
assumptions.  The first one concerns the relation between each set
$\Omega_i$, its covering sets $A_{i,j}$, and its constitutive
constraints $g_{i,j,\ell}$. This assumption states that, if a point is
in the closure of $A_{i,j}$ and satisfies the constitutive constraints
associated with $A_{i,j}$ then this point necessarily belongs to
$\Omega_i$.

\begin{assump} \label{a5}
For all $i \in \{1,\dots,\nblocks\}$ and all $j \in
\{1,\dots,\nops(i)\}$, if $\bx_i \in \overline{A_{i,j}}$ and
$g_{i,j,\ell}(\bx_i) \leq 0$ for $\ell=1,\dots,n_g(i,j)$, then $\bx_i
\in \Omega_i$.
\end{assump}

\begin{figure}[ht!]
\begin{center}
\includegraphics{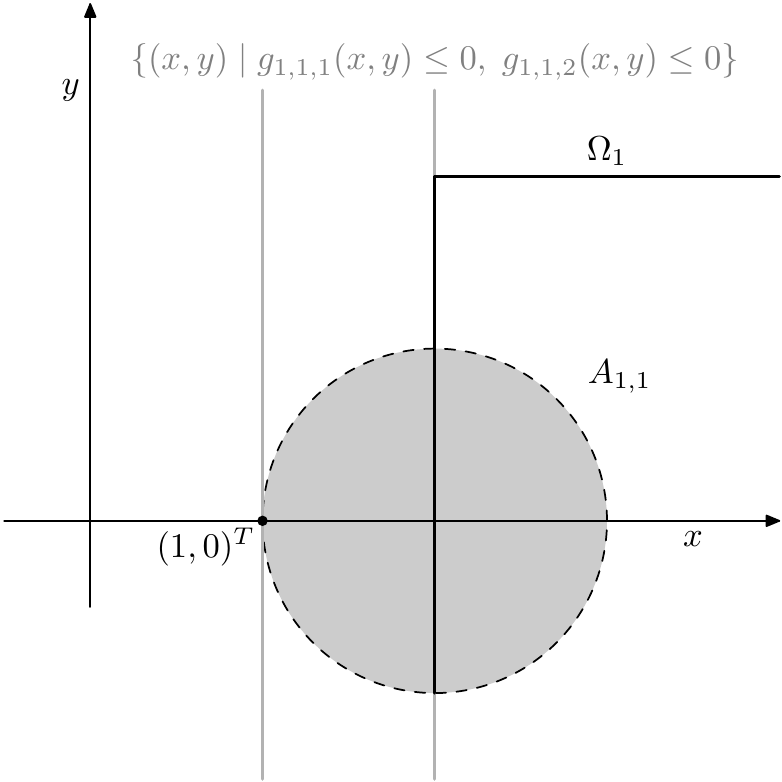}
\end{center}
\caption{The constitutive constraints of $\Omega_1$ over the open set
  $A_{1,1}$ might be given by $g_{1,1,1}(x,y):= (x-1)(x-2) \leq 0$ and
  $g_{1,1,2}(x,y):= - (x-1)(x-2) \leq 0$. This is an example that
  satisfies~(\ref{ineq2}) with $i=j=1$ and $n_g(1,1)=2$ at which
  Assumption~\ref{a5} does not hold, because the point $(x,y)=(1,0)
  \in \overline{A_{1,1}}$ satisfies the constitutive constraints but it does
  not belong to $\Omega_1$.
  }
\label{fig2}
\end{figure}

Figure~\ref{fig2} illustrates a pathological example at which
Assumption~\ref{a5} does not hold.

The second additional assumption states that every global minimizer of
a smooth function onto $\Omega_i$ belongs to some $A_{i,j}$ satisfying
the KKT conditions with respect to the constitutive constraints
related with $A_{i,j}$.

\begin{assump} \label{a6}
For all $i \in \{1,\dots,\nblocks\}$, if $\bx_i$ is a global minimizer
of a smooth function onto $\Omega_i$, there exists $j \in
\{1,\dots,\nops(i)\}$ such that $\bx_i \in A_{i,j}$ and $\bx_i$
satisfies the KKT conditions with respect to the constitutive
constraints associated with $A_{i,j}$.
\end{assump}

Let us show that Assumption~\ref{a6} does not hold in the example of
Figure~\ref{fig1}. Assume that the global minimizer~$G$ of a smooth
function $\psi$ over the set $\Omega_1$ belongs to the boundary of
$\Omega_1$ and to the ball $A_{1,4}$ but does not belong to $A_{1,9}
\cup A_{1,10} \cup A_{1,11} \cup A_{1,12}$ and is not the center
$C_{1,4}$ of $A_{1,4}$\footnote{We already shown, in
  Section~\ref{problem}, that, althought non-regular, if a smooth
  function $\psi$ has a minimizer at $C_{1,4}$, its gradient $\nabla
  \psi(C_{1,4})$ is necessarily null; so that $C_{1,4}$ is a KKT point
  of the minimization of $\psi$ subject to $\varphi \leq 0$.}. For
example, take an adequate infeasible point~$Q$ very close to the
desired global minimizer~$G$ and define $\psi(P) = \|P - Q\|^2$.  The
global minimizer~$G$ belongs only to the open set $A_{1,4}$, the
gradient $\nabla \psi(G)$ is nonnull but, according to the definition
of the constitutive constraint $\varphi$ of $A_{1,4}$, $\nabla
\varphi(G)=0$. Therefore, the KKT condition does not hold in this
case. Fortunately, there exist several simple ways to fix this
drawback. For example, we may define $A_{1,13}$ as a ball with center
in the boundary of $\Omega_1$ such that $C_{1,4}$, the center of
$A_{1,4}$, is on the boundary of this ball and the radius is large
enough so that $A_{1,9} \cap A_{1,13}$ is nonempty. (Analogously, we
may define $A_{1,14}$ as a ball with center in the boundary of
$\Omega_1$ such that $C_{1,4}$ is on the boundary of this ball and the
radius is large enough so that $A_{1,10} \cap A_{1,14}$ is nonempty.)
So, the global minimizer $G$ defined above would belong, not only to
the problematic open set $A_{1,4}$ but also to the newly defined
$A_{1,13}$ (or $A_{1,14}$) where only one linear constitutive
constraint is present and, consequently, KKT necessarily holds.

In order to pursue Alternative 1, at Step~2 of
Algorithm~\ref{bcsmethod}.1, $j_k \in \{1,\dots,\nops(i_k)\}$ and
$\bx_{i_k}^{\trial} \in A_{i_k,j_k}$ must be found such
that~(\ref{bajar}) holds and such that there exist $\mu_{i_k,j_k,\ell}
\geq 0$ for $\ell = 1,\dots, \nge(i_k,j_k)$ for which~(\ref{grachico})
and~(\ref{lambdamu}) hold. To accomplish this, for $j$ from $1$ to
$\nops(i_k)$, provided it is affordable, we could compute a global
minimizer $\boldsymbol{z}_j^*$ of
\begin{equation} \label{qprob}
\begin{array}{c}
\displaystyle
\Minimize_{\bx \in \R^{n_{i_k}}} \nabla_{i_k} f(x^k)^T (\bx -
\bx_{i_k}^k) + \half (\bx - \bx_{i_k}^k) B_k (\bx - \bx_{i_k}^k) +
\frac{\sigma}{2} \| \bx - \bx_{i_k}^k \|^2\\[2mm]
\mbox{subject to } \bx \in \overline{A_{i_k,j}} \mbox{ and } g_{i_k,j,\ell}(\bx)
\leq 0 \mbox{ for } \ell=1,\dots,n_g(i_k,j).
\end{array}
\end{equation}
If the objective function value at $\boldsymbol{z}_j^*$ is
non-positive and $\boldsymbol{z}_j^* \in A_{i_k,j}$, then, defining
$j_k=j$, we have that $\bx_{i_k}^{\trial} = \boldsymbol{z}_j^*$
satisfies~(\ref{bajar}). If we additionaly assume that global
minimizers of~(\ref{qprob}) for every $j \in \{1,\dots,\nops(i_k)\}$
satisfy KKT conditions, then we have that there exist
$\mu_{i_k,j_k,\ell} \geq 0$ for $\ell = 1,\dots, \nge(i_k,j_k)$ for
which~(\ref{grachico}) and~(\ref{lambdamu}) hold.

If none of the global minimizers $\boldsymbol{z}_j^*$ is such the
objective function value at $\boldsymbol{z}_j^*$ is non-positive and
$\boldsymbol{z}_j^* \in A_{i_k,j}$, then we can define $j_k=j^a$ and
$\bx_{i_k}^{\trial} = \boldsymbol{z}_{j^b}^*$, where
$\boldsymbol{z}_{j^b}^*$ is the global minimizer (among the
$\nops(i_k)$ computed global minimizers $\boldsymbol{z}_1^*, \dots,
\boldsymbol{z}_{\nops(i_k)}^*$) that achieves the lowest functional
value of the objective function in~(\ref{qprob}) and $j^a$ is such
that $\boldsymbol{z}_{j^b}^* \in A_{i_k,j^a}$. The functional value of
the objective function of~(\ref{qprob}) at $\boldsymbol{z}_{j^b}^*$ is
non-positive because the objective function vanishes at $x_{i_k}^k$
that is a feasible point of~(\ref{qprob}) for at least one $j \in
\{1,\dots,\nops(i_k)\}$; and $\boldsymbol{z}_{j^b}^* \in A_{i_k,j^a}$
for some $j^a$ because, by Assumption~\ref{a5},
$\boldsymbol{z}_{j^b}^* \in \Omega_{i_k}$ and, by definition,
$\Omega_{i_k} \subset A_{i_k} = \cup_{j=1}^{\nops(i_k)}
A_{i_k,j}$. Moreover, $\boldsymbol{z}_{j^b}^*$ must also be a global
minimizer of~(\ref{qprob}) with $j=j^a$. Thus, by Assumption~\ref{a6},
it fulfills KKT conditions and, therefore, there exist
$\mu_{i_k,j_k,\ell} \geq 0$ for $\ell = 1,\dots, \nge(i_k,j_k)$ for
which~(\ref{grachico}) and~(\ref{lambdamu}) hold.

It is worth noting that the objective function in~(\ref{qprob}) is a
linear function if $B_k=0$ and $\sigma=0$; and it is a convex
quadratic function if $B_k + \sigma I$ is positive definite. Moreover,
it is always possible to choose the open covering sets $A_{i,j}$ for
all $i$ and $j$ in such a way their closures are simple sets like
balls, boxes, or polyhedrons. Furthermore, it is also possible that
more efficient problem-dependent alternatives exist for the
computation of $\bx_{i_k}^{\trial}$; and it is also possible trying to
find $\bx_{i_k}^{\trial}$ satisfying Alternative~2 when $\sigma=0$.

\section{Experiments} \label{experiments}

In this section we describe numerical experiments using the BCD
method. Section~\ref{cTSP} describes what we have called the
continuous version of the traveling salesman problem (TSP). In this
problem, the BCD method is used to evaluate the merit of the function
that should be minimized. Since this function is computed many times
along the whole process, this problem provides many experimental
applications of BCD. In Section~\ref{discropt}, we describe a simple
heuristic and a way to generate a starting point for solving the
continuous TSP. Although simple, these considered methods are part of
the state of the art of methods used to solve the classical
TSP. Moreover, they serve to illustrate the application of the BCD
method, which could be used in the same way in combination with any
other strategy. Section~\ref{heuri} describes a problem-dependent way
to find a~$x_{i_k}^{\trial}$ in the BCD method that satisfies the
requirements of Alternative~2. Section~\ref{numexp} describes the
computational experiment itself.

\subsection{Continuous traveling salesman problem} \label{cTSP}

The travelling salesman problem (TSP) is one of the most studied
combinatorial optimization problems for which a vast literature
exists; see, for example, \cite{applegate}, and the references
therein. In its classical version, $p$ cities with known pairwise
``distances'' $d_{ij}>0$ are given and the problem consists in finding
a permutation $i_1, i_2, \dots, i_p$ that minimizes $d_{i_p,i_1} +
\sum_{\nu=1}^{p-1} d_{i_\nu,i_{\nu+1}}$. In the present work, we
consider a continuous variant of the classical TSP in which ``cities''
are not fixed and, therefore, their pairwise distances vary. More
precisely, given a set of polygons $\Omega_1, \Omega_2, \dots,
\Omega_p$, that may be nonconvex, the problem consists of finding
points $\bx_i \in \Omega_i$ for $i=1,\dots,p$ and a permutation $i_1,
i_2, \dots, i_p$ that minimize $\| \bx_{i_p} - \bx_{i_1} \| +
\sum_{\nu=1}^{p-1} \| \bx_{i_\nu} - \bx_{i_{\nu+1}} \|$. Polygons may
be seen as representing countries, regions, districts, or
neighbourhoods of a city; and the interpretation is that ``visiting a
polygon'' is equivalent to ``visiting any point within the polygon''.

The application of the BCD method in this context is very natural. Any
method to solve the classical TSP requires to evaluate the merit of a
permutation $i_1, i_2, \dots, i_p$ by calculating $d_{i_p,i_1} +
\sum_{\nu=1}^{p-1} d_{i_\nu,i_{\nu+1}}$. In the variant we are
considering, given a permutation $i_1, i_2, \dots, i_p$, the BCD
method is used to find the $\bx_{i_\nu} \in \Omega_{i_\nu}$ for
$\nu=1,\dots,p$ that minimize $\| \bx_{i_p} - \bx_{i_1} \| +
\sum_{\nu=1}^{p-1} \| \bx_{i_\nu} - \bx_{i_{\nu+1}} \|$. In other
words, given a permutation $i_1, i_2, \dots, i_p$, the BCD method is
used to find a solution $x^*$ to the problem
\begin{equation} \label{BCDTSPproblem}
\Minimize_{x \in \R^n} f(i_1,\dots,i_p;x) := \| \bx_{i_p} - \bx_{i_1}
\| + \sum_{\nu=1}^{p-1} \| \bx_{i_\nu} - \bx_{i_{\nu+1}} \| \mbox{
  subject to } \bx_{i_\nu} \in \Omega_{i_\nu} \mbox{ for }
\nu=1,\dots,p,
\end{equation}
where $\nblocks=p$, $n_i=2$ for $i=1,\dots,\nblocks$, $n=2p$, and $x =
(\bx_1^T,\dots,\bx_p^T)^T$; while the problem as a whole consists in
finding a permutation $i_1^*,\dots,i_p^*$ such that
$f(i_1^*,\dots,i_p^*; x^*)$ is as small as possible. That is, the BCD
integrates the process of evaluating the merit of a given
permutation. With this tool, constructive heuristics and
neighborhood-based local searches already developed for the classical
TSP can be adapted to the problem under consideration.

\subsection{Discrete optimization strategy} \label{discropt}

In the present work, among the huge range of possibilities and in
order to illustrate the usage of the BCD method, we consider a local
search with an insertion-based neighborhood. The initial solution is
given by a constructive heuristic also based on insertions, as we now
describe; see~\cite{aarts} and the references therein. The
construction of the initial guess starts defining $(i_1,i_2)=(1,2)$
and $\bx_{i_1} \in \Omega_{i_1}$ and $\bx_{i_2} \in \Omega_{i_2}$ as
the ones that minimize $\| \bx_{i_1} - \bx_{i_2} \|$, computed with
the BCD method. Then, to construct $(i_1,i_2,i_3)$, the method
considers inserting index~$3$ before~$i_1$, between~$i_1$ and~$i_2$,
and after~$i_2$. For each of the three possibilities, optimal
$\bx_{i_1} \in \Omega_{i_1}$, $\bx_{i_2} \in \Omega_{i_2}$, and
$\bx_{i_3} \in \Omega_{i_3}$ are computed with the BCD method. Among
the three permutations, the one with smallest $\| \bx_{i_1} -
\bx_{i_2} \| + \| \bx_{i_2} - \bx_{i_3} \| + \| \bx_{i_3} - \bx_{i_1}
\|$ is chosen. The method proceeds in this way until a permutation
with $p$ elements, that constitutes the initial guess, is completed. A
typical iteration of the local search proceeds as follows. Given the
current permutation $(i_1, i_2, \dots, i_p)$ and its associated points
$\bx_{i_\nu} \in \Omega_{i_\nu}$ for $\nu=1,\dots,p$, each $i_s$ for
$s=1,\dots,p$ is removed and reinserted at all possible places $t \neq
s$. For each possible insertion, corresponding $\bx_{i_1}, \bx_{i_2},
\dots, \bx_p$ are computed with the BCD method. This type of movement
is also known as \textit{relocation} and, as mentioned
in~\cite[p.342]{aarts}, it has been used with great success in the
TSP~\cite{laporte}. Once an insertion is found that improves the
current solution, the iteration is completed, i.e.\ the first
neighbour that improves the current solution defines the new iterate,
in constrast to a ``best movement'' strategy in which all neighbors
are considered and the best of them defines the new iterate. The local
search ends when no neighbour is found that improves the current
iterate.

\subsection{Finding optimal points with BCD method for a given permutation} \label{heuri}

In this section we describe how to solve problem~(\ref{BCDTSPproblem})
with the BCD method. At iteration~$k$ of the BCD method, an index $i_k
\in \{1,\dots,\nblocks\}$ is chosen at Step~1. Then at Step~2, there
are two alternatives. If $\sigma=0$, then $\bx_{i_k}^{\trial}$
satisfying Alternative~1: (\ref{bajar},\ref{grachico},\ref{lambdamu})
or Alternative~2: (\ref{grachi2},\ref{lamb2}) must be computed; while,
if $\sigma>0$, then only Alternative~1 is a
possibility. Section~\ref{subpros2} describes a way of computing
$\bx_{i_k}^{\trial}$ satisfying Alternative~1 for any value of
$\sigma$. However, for the particular problem under consideration,
minimizing $f(x)$ as a function of $\bx_{i_k} \in \Omega_{i_k}$
reduces to
\begin{equation} \label{problemafacil}
\Minimize_{\bx_{i_k} \in \R^2} \| a - \bx_{i_k} \| + \| \bx_{i_k} - b \|
\mbox{ subject to } \bx_{i_k} \in \Omega_{i_k},
\end{equation}
where $a$ and $b \in \R^2$ stand for the ``previous'' and the ``next''
point in the permutation; that, in general, correspond to
$\bx_{i_{k-1}}$ and $\bx_{i_{k+1}}$, respectively. Thus, when
$\sigma=0$, it is easy, computationally tractable, and affordable to
compute the global minimizer of~(\ref{problemafacil}), which clearly
satisfies the requirements of Alternative~2. The global minimizer is
either on the segment $[a,b]$ intersected with $\Omega_{i_k}$ (that
intersection is given by a finite set of segments) or on the boundary
of $\Omega_{i_k}$, which is also given by a finite set of segments
(its edges). Each segment can be parameterized with a single variable
$\lambda \in [0,1]$. Then, the global minimizer
of~(\ref{problemafacil}) is given by the best global minimizer among
the global minimizers of these simple box-constrained one-dimensional
problems. The global minimizer of each box-constrained one-dimensional
problem can be computed with brute force up any desired
precision. Moreover, if multiple solutions exist, in order in increase
the chance of satisfying~(\ref{armijo}), the closest one to
$x_{i_k}^k$ should be preferred.

\subsection{Traveling in S\~ao Paulo City} \label{numexp}

For the numerical experiments, we implemented the discrete
optimization strategy described in Section~\ref{discropt} and the BCD
method (Algorithm~\ref{bcsmethod}.1) described in
Section~\ref{bcsmethod} with the strategy described in
Section~\ref{heuri} for the computation of~$\bx_{i_k}^{\trial}$. In
Algorithm~\ref{bcsmethod}.1, we chose $i_k=mod(k+1,\nblocks)$ and,
based on the theoretical results, we stop the method at iteration $k$,
if $x^k = x^{k-1} = \dots = x^{k-\nblocks+1}$. In the numerical
experiments, following~\cite{aabmm,bmbunch,bmnewgen}, we consider
$\alpha=10^{-8}$. In all our experiments the required conditions at
Step~2 were satisfied using Alternative~2.

All methods were implemented in Fortran~90. Tests were conducted on a
computer with a 3.4 GHz Intel Core i5 processor and 8GB 1600 MHz DDR3
RAM memory, running macOS Mojave (version 10.14.6). Code was compiled
by the GFortran compiler of GCC (version 8.2.0) with the -O3
optimization directive enabled.

The city of S\~ao Paulo, with more than 15 million square kilometers
of extension and more than 12 million inhabitants, is the most
populous city in Brazil, the American continent, the
Portuguese-speaking countries and the entire southern hemisphere. It
is administratively divided into thirty-two regions, each of which, in
turn, is divided into districts, the latter sometimes subdivided into
subdistricts (popularly called neighborhoods); see
\url{https://pt.wikipedia.org/wiki/S%C3%A3o_Paulo}. The city has a
  total of 96 neighborhoods. The considered problem consists in
  finding a shortest route to visit all of them.

\begin{figure}[ht!]
\begin{center}
\includegraphics[scale=0.5]{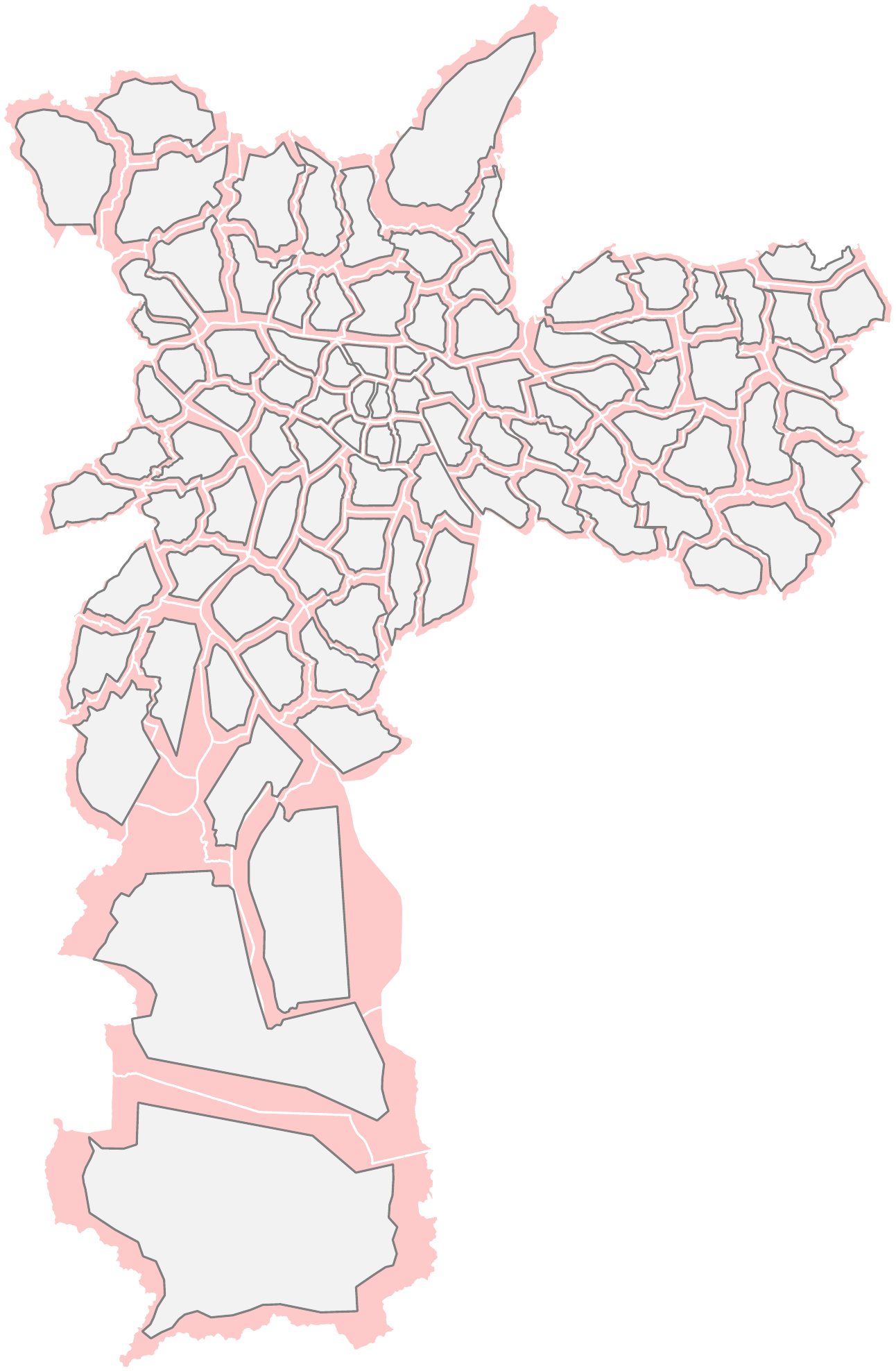}
\end{center}
\caption{Representation of the $p=96$ polygons that determine
  problem~(\ref{BCDTSPproblem}). The considered polygons appear in
  gray, while the original polygons appear in coral on the background
  merely to improve the artistic appearance of the drawing.}
\label{fig3}
\end{figure}

The construction of the problem started by downloading a
political-administrative map of the city from the city hall website;
see \url{http://geosampa.prefeitura.sp.gov.br/}. The map describes
each neighborhood as a polygon. The polygon with more vertices has
5{,}691 vertices and, all together, the polygons have 156{,}852
vertices. To turn the problem into something more tractable, we
redefine the polygons with number of vertices $n_v>100$ (all of them
in fact) by considering only the vertices with indices of the form
form $1 + \lfloor 50 / n_v \rfloor j$ for $j=0,1,2,\dots$. This way,
all polygons were left with a number of vertices between 51 and 57,
totaling 4{,}966 vertices. Moreover, for artistic reasons related to
the graphical representation of the problem, we shrunk each polygon by
20\%. The shrinkage consisted in replacing each vertex $v_i$ by $o_i +
0.8 ( v_i - o_i )$, where the offset $o_i = \half ( x_{\min} +
x_{\max}, y_{\min} + y_{\max})^T$, and $(x_{\min},y_{\min})^T$ and
$(x_{\max},y_{\max})^T$ correspond to the lower-left and upper-right
corners of the smallest rectangle that encloses the polygon. With this
procedure we ended up with the $p=96$ polygons $\Omega_i$ for
$i=1,\dots,p$ that determine the problem~(\ref{BCDTSPproblem}) under
consideration; see Figure~\ref{fig3}.

Table~\ref{tab1} shows the details of the optimization process. The
table shows, for each iteration, the length of the current route. It
also shows, for each iteration, how many neighbors had to be evaluated
to find one that improves the current route. Naturally, each
evaluation of a neighbor corresponds to a call to the BCD
method. Therefore, the next two columns show the number of calls to
the BCD method per iteration and the number of cycles these calls
used. The last two columns of the table show these two values
accumulated over the iterations. It can be noted from the table that
the BCD method is used to solve more than 200{,}000 subproblems and
that this requires, altogether, the execution of more than 3 million
cycles, i.e. an average of 15 cycles per problem. The instance under
consideration has $p=96$ points. The constructive heuristic used to
generate the initial point evaluates (by calling the BCD method)
$O(\half p^2)$ permutations; while the reinsertion neighborhood
evaluates, in the worst case, $O(p^2)$ neighbors. The first line of
the table shows the cost of the constructive heuristic, and is
consistent with what we have just mentioned. The remaining lines show
that in the first 4 iterations and in a few intermediate iterations
the method quickly finds a neighbor that improves the current
solution. On the other hand, the average number of neighbors evaluated
per iteration is 4{,}164, which corresponds to approximately 45\% of
the neighbors. The running time of the algorithm is directly
proportional and totally dependent on the cost of computing
$\bx_{i_k}^{trial}$. The constructive heuristic used to compute the
initial point $x^0$ used 29.58 seconds of CPU time; while the method
as a whole consumed practically one hour of CPU time, exactly
3{,}589.31 seconds.
  
\begin{table}[ht!]
\begin{center}
{\scriptsize
\begin{tabular}{|c|c|cc|cc|}
\hline
\multirow{2}{*}{iter} & \multirow{2}{*}{Route lenght} &
\multicolumn{2}{c|}{Usage of BCD method per iter} &
\multicolumn{2}{c|}{Accumulated usage of BCD method} \\
\cline{3-6}
& & \# calls & \# cycles & \# calls & \# cycles \\
\hline
\hline
 0 & 229{,}139.65 &  4{,}653 &    40{,}018 &    4{,}653 &    40{,}018\\
 1 & 227{,}965.59 &        2 &          26 &    4{,}655 &    40{,}044\\
 2 & 227{,}110.10 &        1 &          10 &    4{,}656 &    40{,}054\\
 3 & 226{,}970.07 &        3 &          56 &    4{,}659 &    40{,}110\\
 4 & 226{,}970.07 &        2 &          24 &    4{,}661 &    40{,}134\\
 5 & 226{,}588.08 &  4{,}125 &    60{,}955 &    8{,}786 &   101{,}089\\
 6 & 226{,}586.52 &  4{,}418 &    66{,}454 &   13{,}204 &   167{,}543\\
 7 & 226{,}575.05 &  4{,}220 &    64{,}793 &   17{,}424 &   232{,}336\\
 8 & 226{,}575.05 &       95 &     1{,}655 &   17{,}519 &   233{,}991\\
 9 & 226{,}573.77 &  4{,}510 &    66{,}294 &   22{,}029 &   300{,}285\\
10 & 226{,}573.77 &      189 &     2{,}658 &   22{,}218 &   302{,}943\\
11 & 226{,}391.25 &  4{,}708 &    69{,}757 &   26{,}926 &   372{,}700\\
12 & 226{,}063.02 &  4{,}789 &    70{,}987 &   31{,}715 &   443{,}687\\
13 & 224{,}708.13 &  3{,}653 &    57{,}797 &   35{,}368 &   501{,}484\\
14 & 224{,}391.45 &  4{,}889 &    73{,}541 &   40{,}257 &   575{,}025\\
15 & 224{,}391.45 &  3{,}747 &    59{,}159 &   44{,}004 &   634{,}184\\
16 & 224{,}236.40 &  4{,}801 &    70{,}928 &   48{,}805 &   705{,}112\\
17 & 224{,}128.38 &  4{,}983 &    74{,}453 &   53{,}788 &   779{,}565\\
18 & 224{,}128.38 &       95 &     1{,}761 &   53{,}883 &   781{,}326\\
19 & 224{,}128.38 &       95 &         908 &   53{,}978 &   782{,}234\\
20 & 224{,}100.64 &  3{,}851 &    58{,}075 &   57{,}829 &   840{,}309\\
21 & 224{,}100.64 &  3{,}838 &    57{,}462 &   61{,}667 &   897{,}771\\
22 & 223{,}681.92 &  4{,}887 &    75{,}693 &   66{,}554 &   973{,}464\\
23 & 223{,}261.18 &  3{,}947 &    66{,}182 &   70{,}501 &  1{,}039{,}646\\
24 & 223{,}261.18 &       95 &         847 &   70{,}596 &  1{,}040{,}493\\
25 & 223{,}242.32 &  5{,}076 &    77{,}901 &   75{,}672 &  1{,}118{,}394\\
26 & 223{,}013.14 &  5{,}572 &    92{,}531 &   81{,}244 &  1{,}210{,}925\\
27 & 221{,}526.66 &  5{,}701 &    91{,}366 &   86{,}945 &  1{,}302{,}291\\
28 & 221{,}469.25 &       94 &         900 &   87{,}039 &  1{,}303{,}191\\
29 & 219{,}846.99 &  5{,}606 &    89{,}792 &   92{,}645 &  1{,}392{,}983\\
30 & 219{,}505.59 &        1 &           4 &   92{,}646 &  1{,}392{,}987\\
31 & 219{,}096.08 &  5{,}761 &    92{,}924 &   98{,}407 &  1{,}485{,}911\\
32 & 218{,}652.50 &  6{,}164 &    96{,}849 &  104{,}571 &  1{,}582{,}760\\
33 & 217{,}719.39 &  6{,}146 &    97{,}136 &  110{,}717 &  1{,}679{,}896\\
34 & 216{,}509.74 &  6{,}256 &    98{,}403 &  116{,}973 &  1{,}778{,}299\\
35 & 215{,}372.59 &  6{,}350 &    99{,}083 &  123{,}323 &  1{,}877{,}382\\
36 & 214{,}676.14 &  5{,}289 &    87{,}276 &  128{,}612 &  1{,}964{,}658\\
37 & 214{,}674.57 &  5{,}077 &    82{,}038 &  133{,}689 &  2{,}046{,}696\\
38 & 214{,}674.57 &  1{,}345 &    29{,}796 &  135{,}034 &  2{,}076{,}492\\
39 & 214{,}102.74 &  3{,}937 &    71{,}897 &  138{,}971 &  2{,}148{,}389\\
40 & 214{,}102.74 &       95 &     1{,}033 &  139{,}066 &  2{,}149{,}422\\
41 & 213{,}754.48 &  5{,}190 &    92{,}893 &  144{,}256 &  2{,}242{,}315\\
42 & 213{,}533.35 &  6{,}049 &    98{,}937 &  150{,}305 &  2{,}341{,}252\\
43 & 213{,}290.76 &  6{,}241 &   102{,}045 &  156{,}546 &  2{,}443{,}297\\
44 & 213{,}231.26 &  7{,}195 &   115{,}060 &  163{,}741 &  2{,}558{,}357\\
45 & 213{,}071.95 &  6{,}338 &    99{,}378 &  170{,}079 &  2{,}657{,}735\\
46 & 213{,}070.90 &  6{,}625 &   103{,}442 &  176{,}704 &  2{,}761{,}177\\
47 & 213{,}032.17 &  7{,}297 &   111{,}568 &  184{,}001 &  2{,}872{,}745\\
48 & 212{,}773.09 &  7{,}681 &   116{,}393 &  191{,}682 &  2{,}989{,}138\\
49 & 212{,}499.29 &  7{,}585 &   115{,}190 &  199{,}267 &  3{,}104{,}328\\
50 & 212{,}292.01 &  8{,}161 &   125{,}796 &  207{,}428 &  3{,}230{,}124\\
51 & 212{,}292.01 &  9{,}120 &   154{,}558 &  216{,}548 &  3{,}384{,}682\\
\hline
\end{tabular}}
\end{center}
\caption{Performance of the heuristic method applied to solve the
  considered instance of the continuous version of the TSP problem.}
\label{tab1}
\end{table}

\begin{figure}[ht!]
\begin{center}
\input{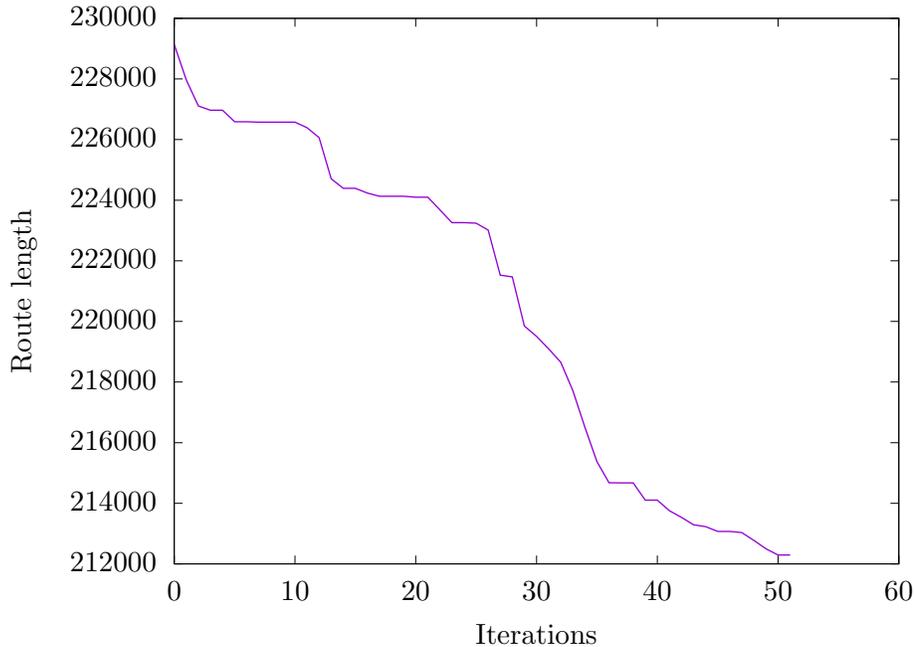}
\end{center}
\caption{Route length as a function of the iteration number.}
\label{fig4}
\end{figure}

\begin{figure}[ht!]
\begin{center}
\begin{tabular}{ccc}
\includegraphics[scale=0.35]{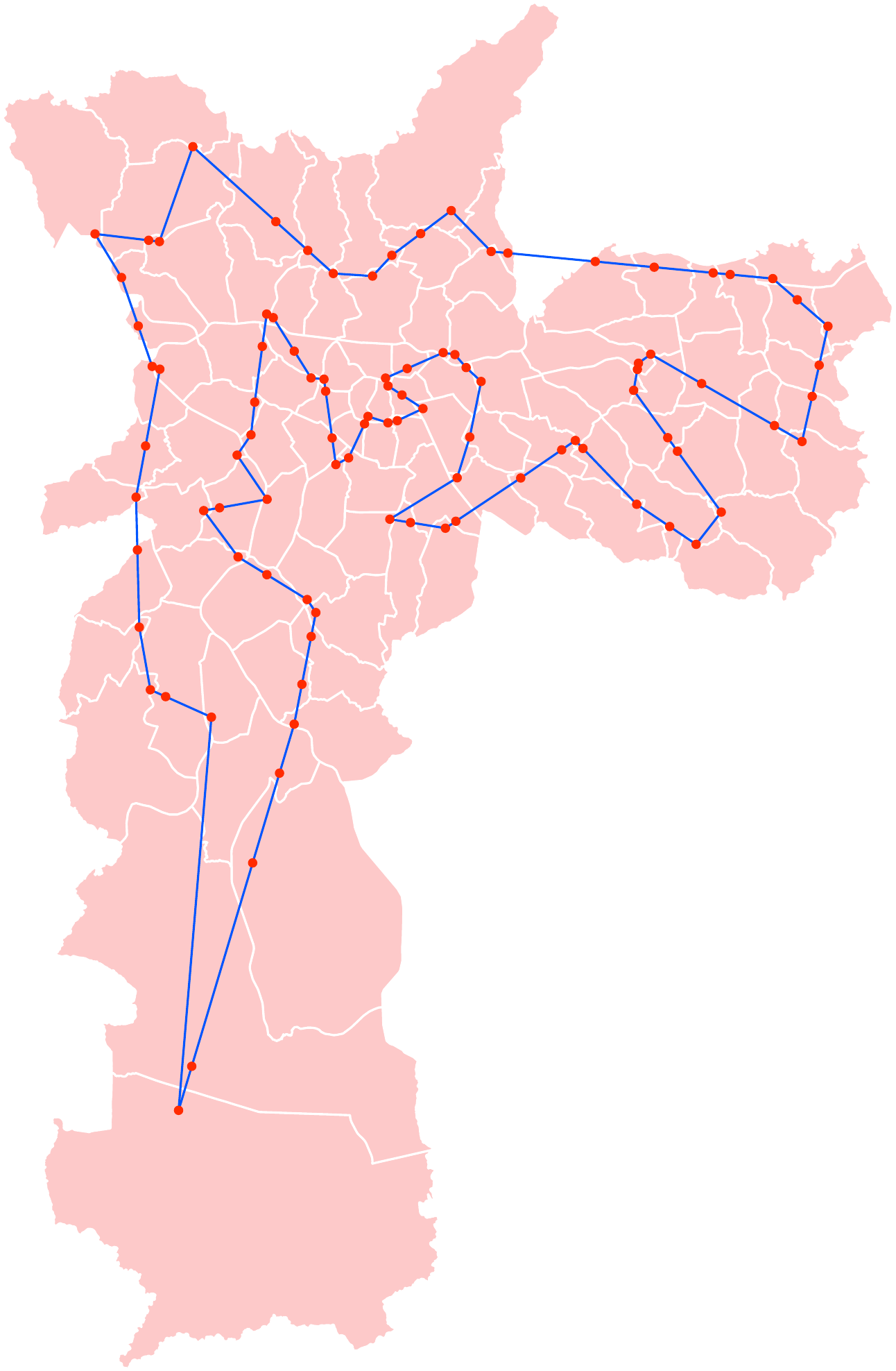} &
\includegraphics[scale=0.35]{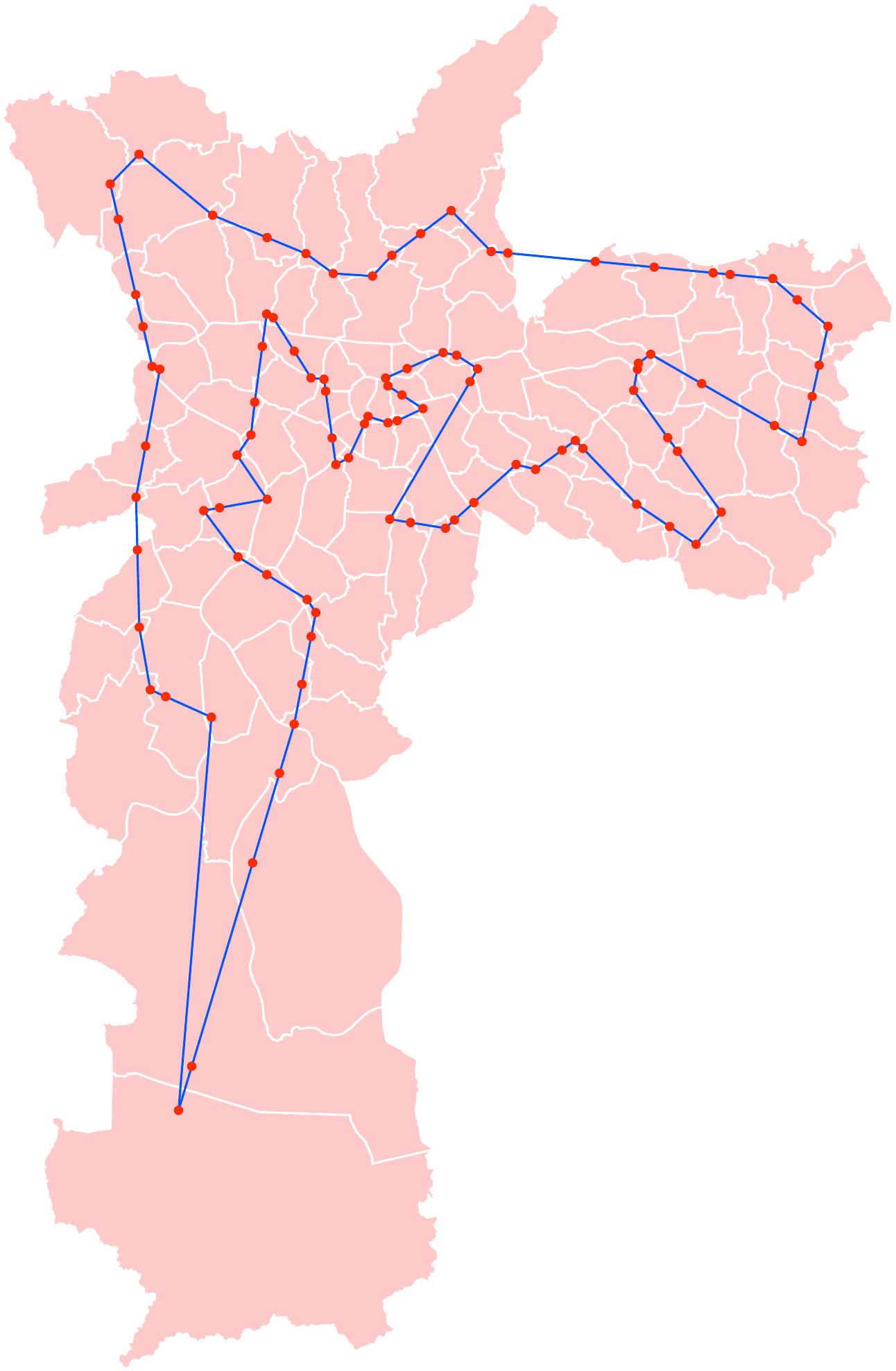} &
\includegraphics[scale=0.35]{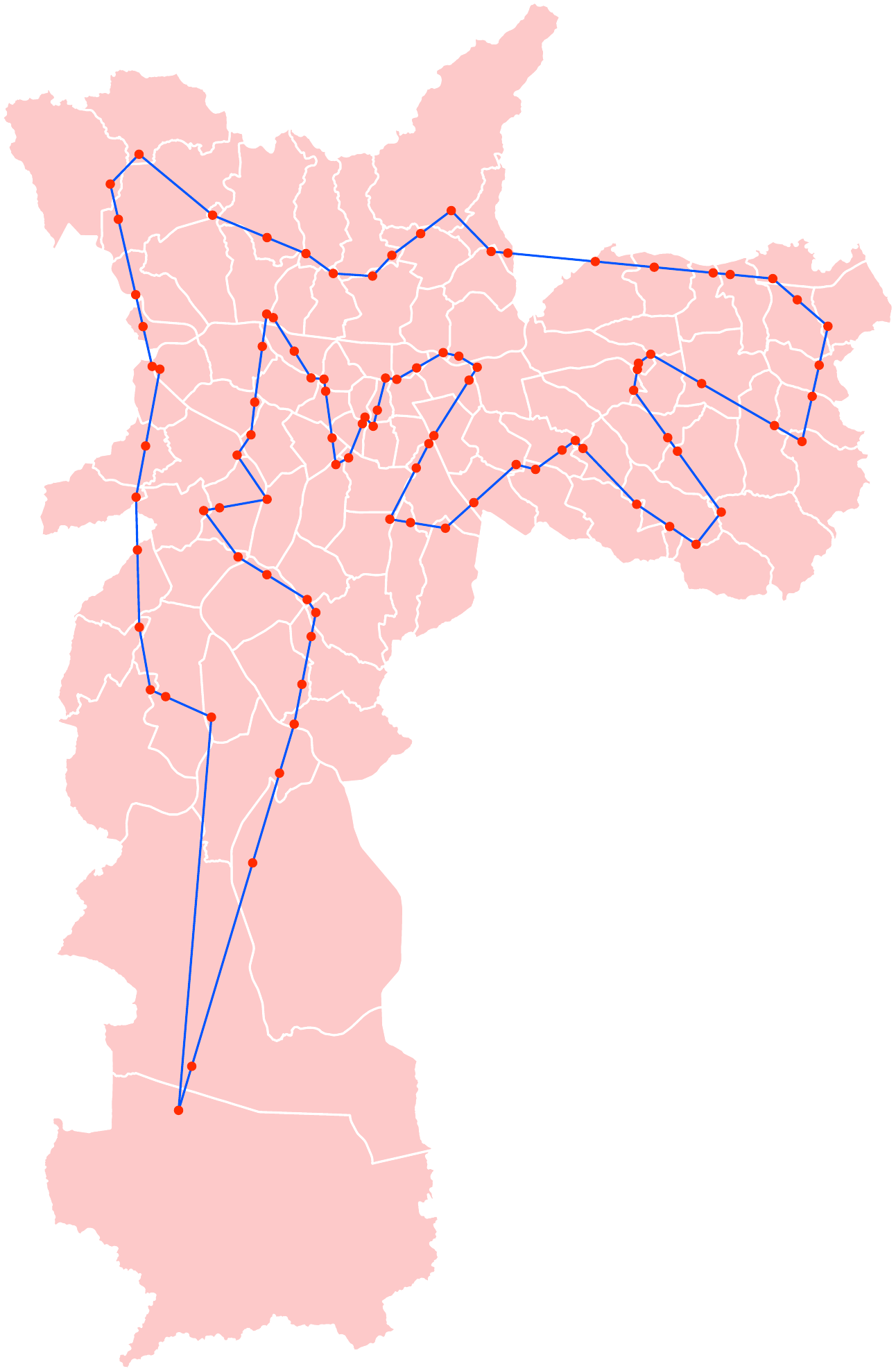} \\
(a) 229{,}139.65 & (b) 226{,}573.767 & (c) 224{,}100.64 \\
\includegraphics[scale=0.35]{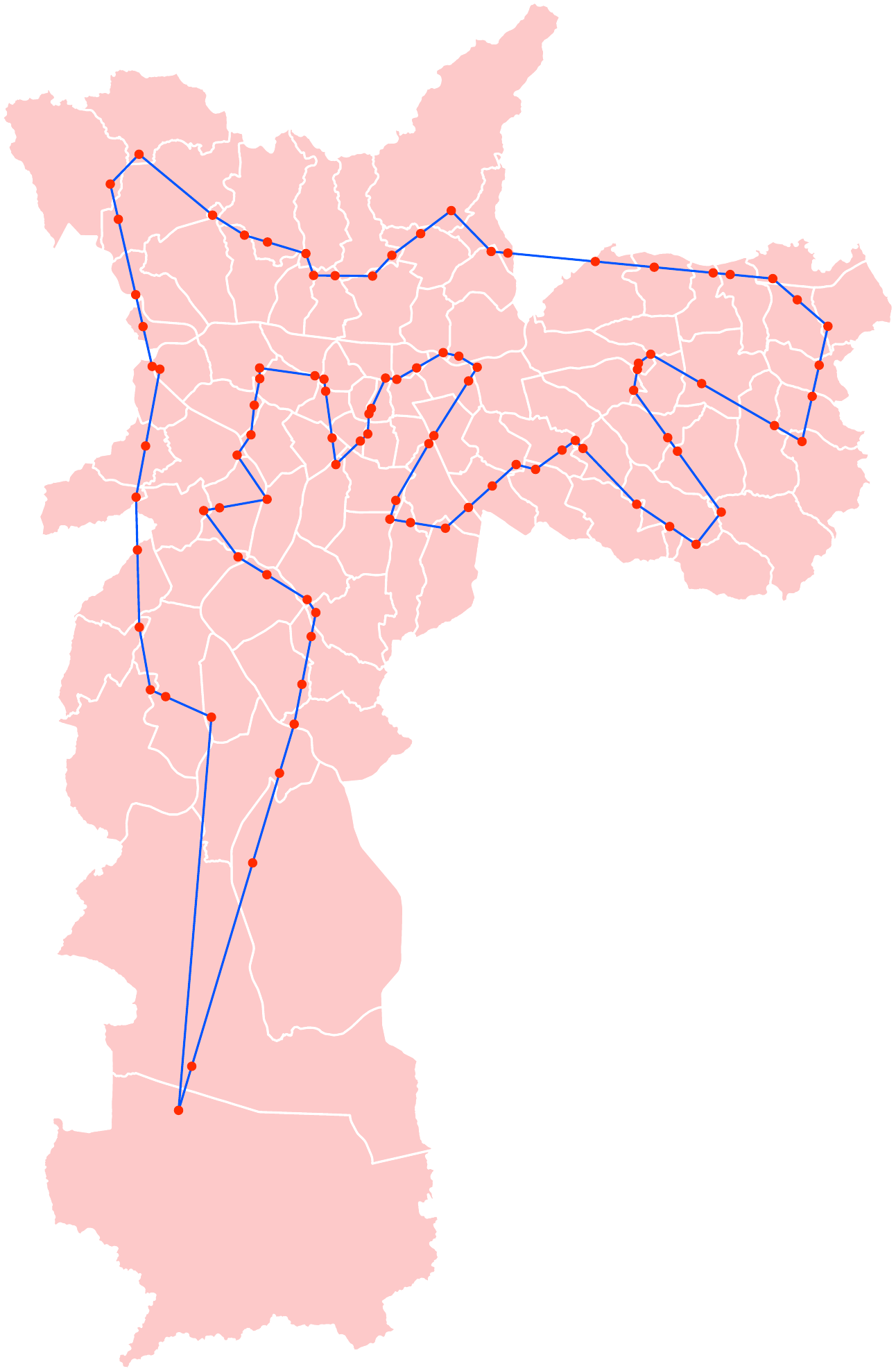} &
\includegraphics[scale=0.35]{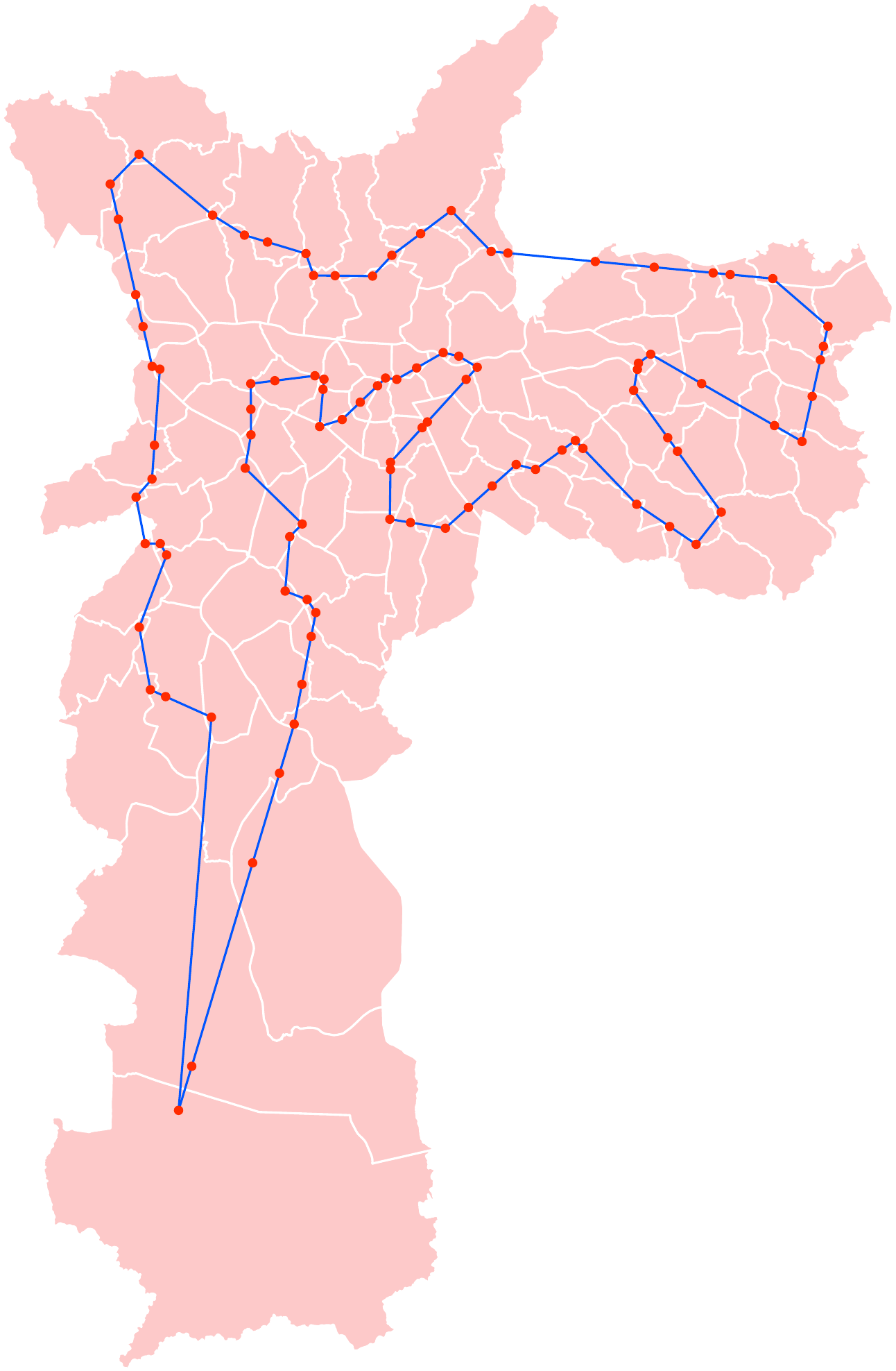} &
\includegraphics[scale=0.35]{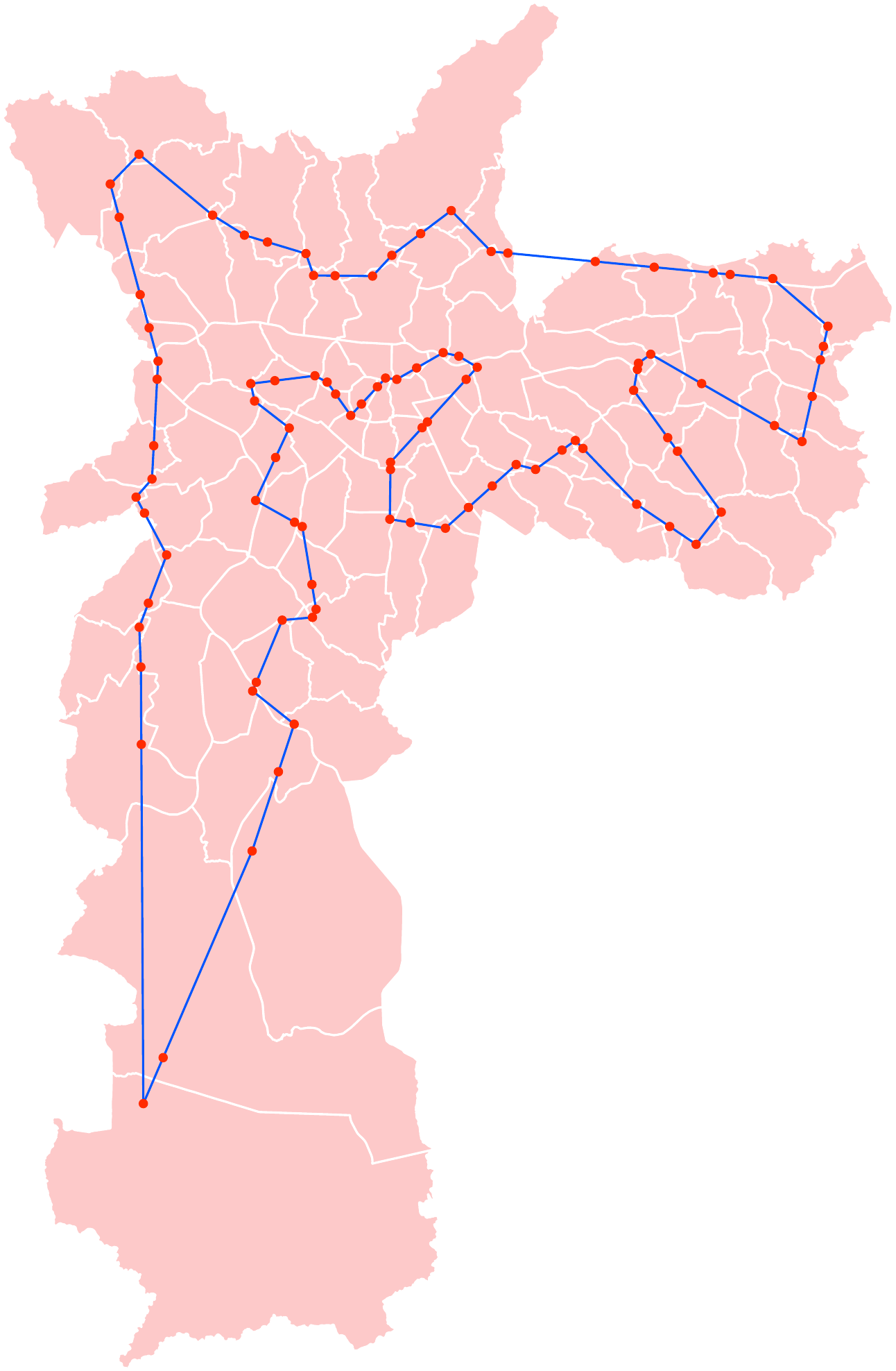} \\
(d) 219{,}505.59 & (e) 214{,}102.74 & (f) 212{,}292.01 \\
\end{tabular}
\end{center}
\caption{Sample of the routes that are built throughout the iterative
  optimization process. The lenght of each route appears near to the
  route. The map of S\~ao Paulo city appears in the background for
  artistic purposes, but the polygons representing the districts are
  being omitted for the sake of clarity. (a) Represents the initial
  guess given by the constructive heuristic; (f) Represents the final
  iterate (obtained at iteration~51); and (b)--(e) Represent the
  iterands of the iterations 10, 20, 30 and 40, respectively. It is
  worth noting that the red dots, each always within its respective
  polygon that is not being displayed, move from one iteration to
  another.}
\label{fig5}
\end{figure}

\begin{figure}[ht!]
\begin{center}
\includegraphics[scale=0.7]{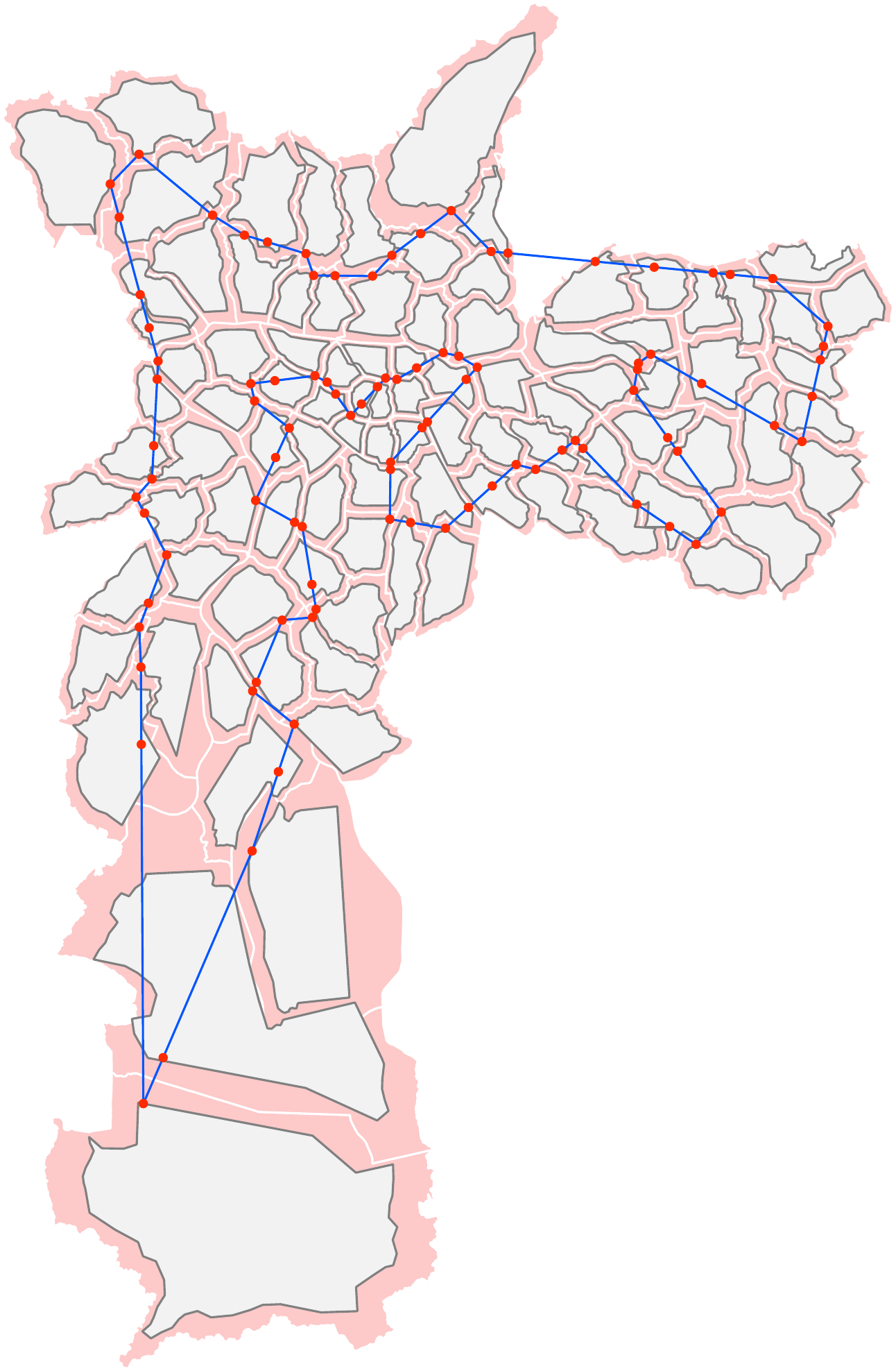}
\end{center}
\caption{Final iterate with route lenght equal to 212{,}292.01.}
\label{fig6}
\end{figure}

Figure~\ref{fig4} shows the evolution of the route length over the
iterations of the method; while Figure~\ref{fig5} shows some of the
generated routes. Figure~\ref{fig6} shows the final iterate in detail.

\begin{figure}[ht!]
\begin{center}
\begin{tabular}{>{\centering\arraybackslash}m{4in} >{\centering\arraybackslash}m{1in}}
\includegraphics[scale=0.7]{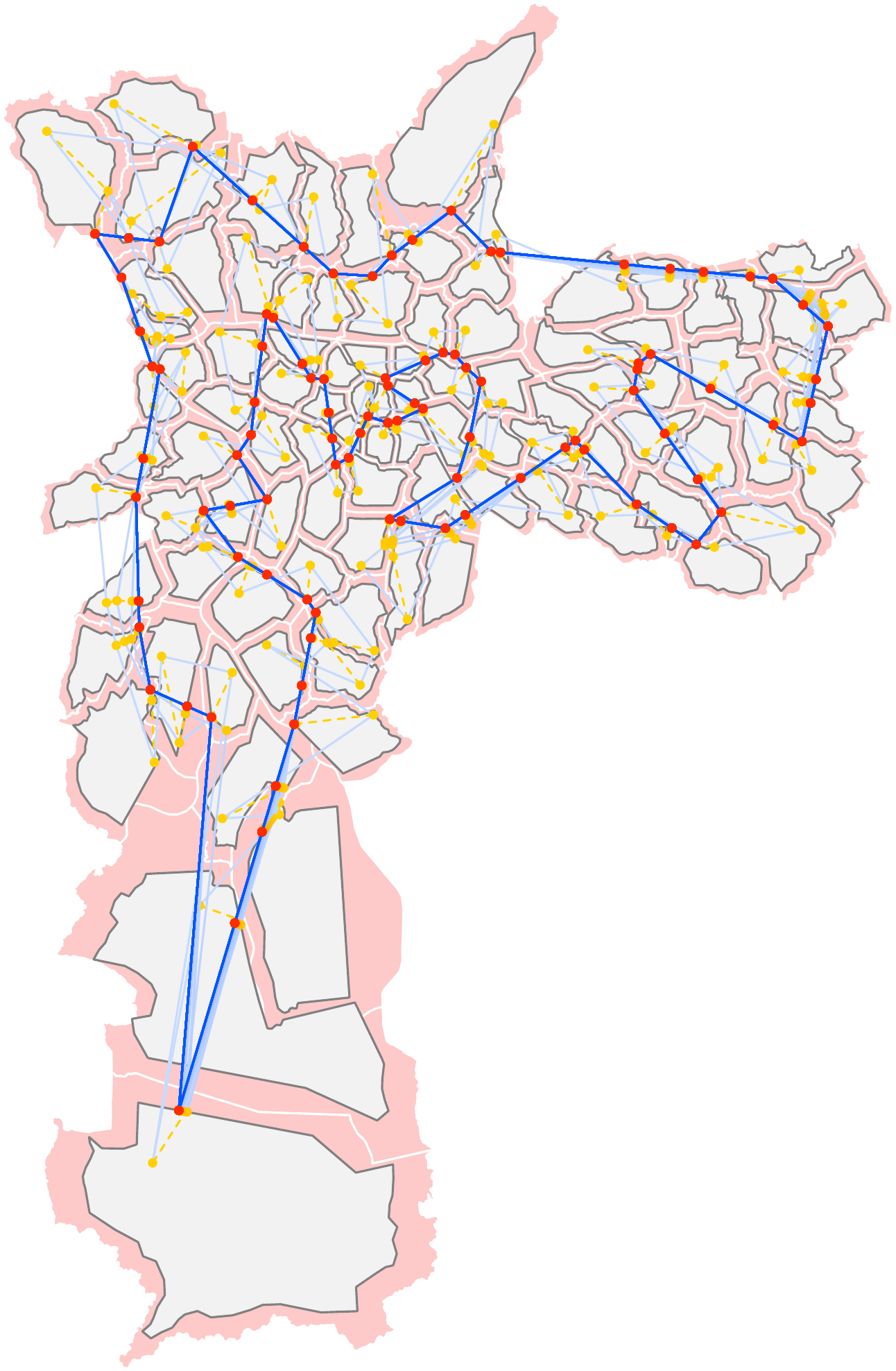} &
\begin{tabular}{|c|c|}
\hline
Cycle & Route lenght \\
\hline
  0 & 425{,}047.86\\
  1 & 244{,}338.54\\
  2 & 232{,}471.73\\
  3 & 230{,}402.57\\
  4 & 229{,}223.98\\
  5 & 229{,}171.84\\
  6 & 229{,}142.28\\
  7 & 229{,}140.84\\
  8 & 229{,}140.22\\
  9 & 229{,}139.93\\
 10 & 229{,}139.79\\
 11 & 229{,}139.72\\
 12 & 229{,}139.68\\
 13 & 229{,}139.66\\
\hline
\end{tabular}
\end{tabular}
\end{center}
\caption{Sequence of iterands of the BCD method when applied to random
  initial points within the polygons and with the order given by the
  constructive heuristic used to generate the initial point.}
\label{fig7}
\end{figure}

Up to this point, we have described and illustrated in detail the way
in which the discrete heuristic method described in
Section~\ref{discropt} made intensive use of the BCD method to solve
the problem~(\ref{BCDTSPproblem}). We close the numerical experiments
section by showing in Figure~\ref{fig7}, with a graphic and a table,
the iterands of the BCD method for a specific fixed permutation. To
make this figure, we considered the permutation given by the
constructive heuristic used to generate $x^0$ and we randomly draw
points inside each of the polygons. The graph and table show the
iterations for 13 complete cycles. The method actually uses 44 cycles,
but the functional value varies from 229{,}139.66 at the end of cycle
13 to 229{,}139.65 at the end of cycle 44, when it stops because all
the variables' blocks are repeated from cycle 43 to cycle 44. The
initial points are in yellow or light orange and the color changes to
red at cycle 13. The evolution of each point is marked with a dotted
line whose color changes along with the color of the
point. Independently of that, the route determined by the points of
each cycle is marked in blue. The route with the lightest blue
corresponds to the route given by the initial points (yellow or light
orange) and the color of the route gets darker and darker until it
reaches the route of cycle 13, marked with the strongest blue. Roughly
speaking, the points move a little more in the first 3 cycles, in
which the objective function decreases the most, and then there are
only small accommodations of the points until the method
converges. The authors are aware of the difficulty to see the figure
clearly; a zoom in the image is recommended to see the details of the
evolution of the iterands. In particular, the middle left part clearly
shows how the route is being modified as the~points~move.

\section{Conclusions} \label{conclusions}

The framework presented in the present work could be extended in order
to consider Taylor-like high-order models \cite{bgmst} satisfying
well-established regularity assumptions, as it has been done in
\cite{aabmm} for the case of box constraints. However, theoretical
results in \cite{aabmm} reveal that using high-order models associated
with Coordinate Descent methods is not worthwhile. The reason is that
overall computed work is dominated by the necessity of obtaining fast
decrease of the distance between consecutive iterates, whereas
high-order models do not help for achieving such purpose.

More interesting, from the practical point of view, is to exploit the
particular case in which the constraints that define each $\Omega_i$
may be expressed in the form of global inequalities and
equalities. (Of course, this is a particular case of the one addressed
in this paper that corresponds to set $\nops(i)=1$ for all
$i=1,\ldots,\nblocks$.)  In this case the obvious choice for solving
subproblems consists of using some well established constrained
optimization software. From the theoretical point of view there is
nothing to be addded, since practical optimization methods for
constrained optimization may fail for different reasons, leading the
abrupt interruption of the overall optimization process. However, we
have no doubts that in many practical problems the standard
constrained optimization approach associated with block coordinate
descent should be useful.

The reason why, in this paper, we considered feasible sets $\Omega_i$
with the local constrained structure defined by open covering sets and
constitutive constraints is not strictly related to block coordinate
methods. In fact, in contact with several practical problems (an
example of which is the one presented in Section~\ref{experiments}) we
observed that the non-global structure of constraints is not unusual
and needs specific ways to be handled properly. We believe that
different approaches than the one suggested in this paper are
possible, most of them motivated by the particular structure of the
practical problems at hand. Further research is expected in the
following years with respect to this subject.

\end{document}